\documentclass{amsart}

\newtheorem{theorem}{Theorem}
\newtheorem{lemma}[theorem]{Lemma}
\newtheorem{proposition}[theorem]{Proposition}

\theoremstyle{definition}

\def\half{{\textstyle{\frac12}}}
\def\third{{\textstyle{\frac13}}}
\def\fourth{{\textstyle{\frac14}}}

\def\threefourths{{\textstyle{\frac34}}}
\def\threehalves{{\textstyle{\frac32}}}
\def\onesixteenth{{\textstyle{\frac1{16}}}}
\def\threesixteenths{{\textstyle{\frac3{16}}}}

\def\pa{\partial}

\def\txtcurl{\text {curl}} 

\def\agamma{ {}^{ \hbox{\tiny $a$} }\!\gamma }
\def\hgamma{ {}^{ \hbox{\tiny $h$} }\!\gamma } 
\def\aG{ {}^{ \hbox{\tiny $a$} }\!G } 
\def\hG{ {}^{ \hbox{\tiny $h$} }\!G } 
\def\aR{ {}^{ \hbox{\tiny $a$} }\!R }
\def\hR{ {}^{ \hbox{\tiny $h$} }\!R } 
\def\aK{ {}^{ \hbox{\tiny $a$} }\!K }
\def\hK{ {}^{ \hbox{\tiny $h$} }\!K } 
\def\ay{ {}^{ \hbox{\tiny $a$} }\!y }

\usepackage{amssymb} 
\title[Zermelo navigation on Riemannian manifolds]
      {Zermelo navigation on Riemannian manifolds} 

\author[Bao, Robles and Shen]
       {DAVID BAO$^*$, COLLEEN ROBLES$^{**}$ and Z. Shen}
\thanks{$^*$Research supported in part by R. Uomini and the  
            S.S. Chern Foundation for Mathematical Research}
\thanks{$^{**}$Research supported in part by the  
               UBC University Graduate Fellowship Program} 



\address{Department of Mathematics, University of Houston \\
         Houston, TX 77204-3008, USA}
\address{Department of Mathematics, University of Rochester \\
         Rochester, NY 14627-0001, USA}
\address{Department of Mathematical Sciences, IUPUI \\ 
         Indianapolis, IN 46202-3216, USA} 

\email{bao@math.uh.edu}
\email{robles@math.rochester.edu}
\email{zshen@math.iupui.edu} 

\begin{document}

%-----------------------------------------------------------------------------

\begin{abstract}
 In this paper, we study Zermelo navigation on Riemannian manifolds and use 
 that to solve a long standing problem in Finsler geometry.  Namely, the 
 complete classification of strongly convex Randers metrics of constant flag 
 curvature. 
\end{abstract}

              %-----------------------------------------------

\maketitle 

%-----------------------------------------------------------------------------

\section{Introduction} 

                %....................................

 \subsection{Purpose}

  We have three goals in this paper.  

  The first is to describe Zermelo's 
  problem of navigation on Riemannian manifolds.  Zermelo aims to find the 
  paths of shortest travel time in a Riemannian manifold $(M,h)$, under the 
  influence of a wind or a current which is represented by a vector field $W$. 
  We point out that the solutions are the geodesics of a strongly convex 
  Finsler metric, which is of Randers type and is {\it necessarily} 
  non-Riemannian unless $W$ is zero.  Conversely, we show constructively that 
  every strongly convex Randers metric arises as the solution to Zermelo's 
  navigational problem on some Riemannian landscape $(M,h)$, under the 
  influence of an appropriate wind $W$.  This is the content of Proposition 1 
  in \S 2.3. 

  Randers metrics are interesting not only as solutions to Zermelo's problem 
  of navigation.  They form a ubiquitous class of metrics with a strong 
  presence in both the theory and applications of Finsler geometry.  Of 
  particular interest are Randers metrics of constant flag curvature, the 
  latter being the Finslerian analog of the Riemannian sectional curvature.  

  It is the second goal of this paper to describe strongly convex Randers 
  metrics of constant flag curvature via Zermelo navigation.  Unlike previous 
  characterization results \cite{BR03, MS02}, the navigation description has 
  the advantage of clearly illuminating the underlying geometry.  More 
  precisely, suppose $(h,W)$ is the navigation data of a strongly convex 
  Randers metric $F$.  Then: $F$ has constant flag curvature $K$ if and only 
  if there exists a constant $\sigma$ such that, $h$ has constant sectional 
  curvature $K + \onesixteenth \sigma^2$ and $W$ satisfies the equation 
  $\mathcal L_W h = - \sigma h$ (namely, $W$ is an infinitesimal homothety of 
  $h$).  This is Theorem 3 in \S 4.4 of the paper. 

  Our third goal has two components. 
  {\bf (1)} Use the navigation description to classify Randers metrics of 
            constant flag curvature.  This problem was proposed by 
            Ingarden about half a century ago.  Until 2002, it was erroneously 
            thought to have been solved by Yasuda--Shimada in 1977; see 
            \cite{BR03, MS02} for references therein.  Our classification 
            lists, explicitly, all the vector fields $W$ that will perturb a 
            Riemannian space form $h$ into a strongly convex Randers metric of 
            constant flag curvature.  This is done in Theorem 7 of \S 6.1. 
  {\bf (2)} Parametrize the moduli space for the equivalence classes of 
            locally isometric $n$-dimensional Randers metrics with constant 
            flag curvature $K$.  In striking contrast with the Riemannian 
            setting -- where the moduli space consists of a single point for 
            each value of $K$ -- the Randers moduli space is of dimension 
            $n/2$ when $n$ is even, and either $(n+1)/2$ or $(n-1)/2$ when $n$ 
            is odd.  The specifics are detailed in Propositions 8 (\S 7.2.5), 
            9 (\S 7.3.1), 10 (\S 7.4.4), respectively for $K$ positive, zero, 
            or negative. 
 
  \medskip 

  Finally, we test the usefulness of the classification by applying it to two 
  special cases.  
  \begin{itemize} 
   \item First, those $W$ which effect projectively flat strongly convex 
         Randers metrics of constant flag curvature $K$ are singled out.  We 
         find that up to local isometry, the {\it non}-Riemannian ones consist 
         of a $1$-parameter family of locally Minkowskian metrics when $K=0$, 
         and a single variant of the Funk metric for each $K < 0$.  In 
         particular, every projectively flat strongly convex Randers metric of 
         constant positive flag curvature must be locally isometric to a 
         standard (Riemannian) sphere.  This discussion constitutes \S 8.3.  
         Our conclusions are also compared with the main result of Shen 
         \cite{S02a}. 
   \item Next, the general classification is specialized to the case in which 
         the tensor $\theta_i := b^s \, \txtcurl_{si}$ vanishes.  This enables 
         us to list explicitly all the Randers metrics addressed by systems of 
         nonlinear partial differential equations in the {\it corrected} 
         Yasuda--Shimada theorem \cite{BR03, MS02}.  Such is the thesis of 
         \S 9.2.  We find that strongly convex {\it non}-Riemannian Randers 
         metrics of constant flag curvature $K$ and $\theta=0$ comprise, 
         {\it up to local isometry}, three small but distinguished camps. 
         \begin{itemize} 
          \item[$\circ$] 
           $K < 0$: there is just a single variant of the Funk metric for 
                    each $K$.  
          \item[$\circ$] 
           $K = 0$: the Yasuda--Shimada theorem tells us that only locally 
                    Minkowski metrics belong to this camp; we show that, in 
                    fact, there is simply a one parameter family of them. 
          \item[$\circ$]
           $K > 0$: this is the most enigmatic case.  There is exactly a 
                    $1$-parameter family of the $\theta=0$ metrics on the odd 
                    dimensional spheres, and none on the even dimensional 
                    spheres.  Furthermore, such conclusion holds whether the 
                    metrics being sought are locally or globally defined.   
         \end{itemize} 
         The classification of the $K > 0$ metrics within the $\theta = 0$ 
         family has previously been done by Bejancu--Farran \cite{BF02, BF03}. 
         However, our description of the local isometry classes offers a 
         perspective which is totally different from theirs.   
  \end{itemize} 
   
                %....................................

 \subsection{Summary of contents}

  Section 2 presents Zermelo's problem of navigation on Riemannian 
  manifolds, and its solution. 

  We specialize to concrete 3-dimensional Riemannian space forms in \S 3.
  These examples are categorized into three subsections, dealing with 
  spheres, Euclidean space, and the Klein model of hyperbolic geometry,
  respectively.  For each model we present examples of Zermelo's 
  navigation which produce Randers metrics of constant flag curvature.
  In \S 3.4, we review the definition of Finsler metrics of constant flag 
  curvature.

  Section 4 begins by recalling a previous characterization result.  It also 
  includes a Matsumoto identity which exhibits the interplay between the 
  constant $\sigma$ (in the equation $\mathcal L_W h = - \sigma h$) and the 
  constant flag curvature $K$.  This is followed by the navigation description 
  of strongly convex Randers metrics of constant flag curvature $K$.  

  Before presenting the classification theorem we pause in \S 5 to derive a 
  complete list of allowable vector fields for each of the three standard 
  models of Riemannian space forms.  With the list in hand, \S 6 gives the 
  classification of strongly convex Randers metrics of constant flag 
  curvature; both local and global aspects are treated.  The isometry classes 
  of constant flag curvature Randers metrics comprise the focus of \S 7.  We 
  make explicit the requisite Lie theory (mostly for a non-compact subgroup of 
  the Lorentz group) in the Appendix (\S 10), and (then) give concrete 
  descriptions of the moduli space and its dimension in \S 7.  Section 8 
  contains a discussion of projectively flat Randers metrics of constant flag 
  curvature.  Finally, in \S 9, we specialize our classification to the 
  $\theta = 0$ case. 

  The dimension counts established in sections 7 through 9 can in essence be 
  summarized by the following table: 
   
   \begin{center} 
   \begin{tabular}{|c|c|c|c|c|c|}
    \hline
    \multicolumn{2}{|c|}{} &
    \multicolumn{4}{|c|}
     {\bf Moduli space's dimension}                                         \\
    \hline
                    & & & & \multicolumn{2}{|c|}{$K < 0$}                   \\ 
                    \cline{5-6}
    \raisebox{1.5ex}[0pt]{CFC metrics} 
                                  & \raisebox{1.5ex}[0pt]{dim $M$} 
                                  & \raisebox{1.5ex}[0pt]{$K > 0$} 
                                  & \raisebox{1.5ex}[0pt]{$K = 0$} 
                                  & $\sigma = 0$ & $\sigma \not= 0$         \\
    \hline
     Riemannian     &                      & \multicolumn{3}{|c|}{} &       \\ 
     $b$ equiv. $W =0$  
                    & \raisebox{1.5ex}[0pt]{$n \geqslant 2$} 
                    & \multicolumn{3}{|c|}{\raisebox{1.5ex}[0pt]{0}}
                    & \raisebox{1.5ex}[0pt]{empty}                          \\ 
    \hline
     Projectively   & & & & & 
                       %\multicolumn{3}{|c|}{} 
                                                                            \\ 
     flat           & \raisebox{1.5ex}[0pt]{$n \geqslant 2$} 
                    & \raisebox{1.5ex}[0pt]{$0^\ast$}
                    & \raisebox{1.5ex}[0pt]{1}
                    & \raisebox{1.5ex}[0pt]{$0^\ast$}
                    & \raisebox{1.5ex}[0pt]{$0^\dagger$}                    \\ 
    \hline
     Yasuda--Shimada & even $n$ & $0^\ast$ & & & \\ \cline{2-3}
     $\theta=0$     & odd $n$  & 1        & \raisebox{1.5ex}[0pt]{1}
                                          & \raisebox{1.5ex}[0pt]
                                            {$0^\ast$}    
                                          & \raisebox{1.5ex}[0pt]
                                            {$0^\dagger$}                   \\ 
    \hline
     Unrestricted   & even $n$ & \multicolumn{4}{|c|}{$n/2$} \\ \cline{2-6}
     Randers        & odd $n$  & \multicolumn{3}{|c|}{$(n+1)/2$}
                               & $(n-1)/2$                                  \\
    \hline
     \multicolumn{6}{|l|}{\it The moduli spaces of dimension 0 
                              consist of a single point.}                   \\
     \multicolumn{6}{|l|}{${}^\ast$ 
                          {The single isometry class is 
                           Riemannian.}}                                    \\
     \multicolumn{6}{|l|}{${}^\dagger$
                          {The single isometry class is 
                           non-Riemannian, of Funk type.}}                  \\ 
    \hline
  \end{tabular}
  \end{center} 
 
                %-----------------------------------

 \subsection{Acknowledgments} 

  We thank A. Bejancu for informing us of his joint results with H. Farran 
  about constant flag curvature Randers metrics with $\theta = 0$ on the 
  $n$-sphere.  We also thank R. Bryant for encouraging us to count isometry 
  classes properly, and for bringing \cite{Br02} to our attention.  The 
  latter provides a framework with which one may hope to classify constant 
  flag curvature Finsler metrics that are not necessarily of Randers type. 

%-----------------------------------------------------------------------------

\section{Zermelo navigation} 

                %------------------------------------

 \subsection{Perturbing Riemannian metrics by vector fields} 

                %------------------------------------

  \subsubsection{Background Riemannian metric and perturbing vector field} 

   Given any Riemannian metric $h$ on a differentiable manifold $M$,  
   denote the corresponding norm-squared of tangent vectors $y \in T_xM$ by 
   \begin{displaymath} 
    |y|^2  
    := h_{ij} \, y^i y^j  
     = h(y,y) \, .   
   \end{displaymath} 
   Think of $|y|$ as the {\it time} it takes, using an engine with 
   a fixed power output, to travel from the base(point) of the vector $y$ 
   to its tip. Note the symmetry property $|-y| = |y|$. 

   The unit tangent sphere in each $T_xM$ consists of all those tangent 
   vectors $u$ such that $|u| = 1$. Now introduce a vector field $W$ such that 
   $|W| < 1$, thought of as the spatial velocity vector of a mild wind on the 
   Riemannian landscape $(M, h)$. Before $W$ sets in, a journey from the base 
   to the tip of any $u$ would take 1 unit of time, say, 1 second. The effect 
   of the wind is to cause the journey to veer off course (or merely off 
   target if $u$ is collinear with $W$). Within the same 1 second, we traverse 
   not $u$ but the resultant $v = u + W$ instead. 

   As an example, suppose $|W| = \half$. If $u$ points along $W$ (that is, 
   $u = 2W$), then $v = \threehalves u$. Alternatively, if $u$ points 
   opposite to $W$ (namely, $u = -2W$), then $v = \half u$. In these two 
   scenarios, $|v|$ equals $\threehalves$ and $\half$ instead of $1$. So, 
   with the wind present, our Riemannian metric $h$ no longer gives the travel 
   time along vectors. This prompts the introduction of a function $F$ on the 
   tangent bundle $TM$, in order to keep track of the travel time needed 
   to traverse tangent vectors $y$ under windy conditions. For all those 
   resultants $v = u + W$ mentioned above, we have $F(v) = 1$. In other words, 
   within each tangent space $T_xM$, the unit sphere of $F$ is simply the 
   $W$-translate of the unit sphere of $h$.  Since this $W$-translate is no 
   longer centrally symmetric, $F$ cannot possibly be Riemannian. 

               %....................................

  \subsubsection{Formula for the new norm $F$} 

   Start with the fact $|u| = 1$; equivalently, $h(u,u) = 1$. Into this, 
   we substitute $u = v - W$ and then 
   $h(v,W) = |v| \, |W| \cos \theta$. After using the abbreviation 
   $\lambda := 1 - |W|^2$ to reduce clutter, we have
   $|v|^2 - \left ( 2 \, |W| \cos \theta \right ) |v| - \lambda = 0$. 
   Since $|W| < 1$, the resultant $v$ is never zero, hence $|v| > 0$. 
   This leads to  
   $|v| = |W| \cos \theta + \sqrt{ \, |W|^2 \cos^2 \theta + \lambda \, }$, 
   which we abbreviate as $p + q$. Since $F(v) = 1$, 
   we see that 
   \begin{displaymath} 
    F(v) = 1                                                       
         = |v| \frac{1}{ \ q + p \ } 
         = |v| \frac{q - p}{ \ q^2 - p^2 \ }                         
         =  \frac{ \sqrt{ \, [h(W,v)]^2 + |v|^2 \lambda \, } }{\lambda} 
          - \frac{ \ h(W,v) \ }{\lambda} \, . 
   \end{displaymath} 
   
   It remains to deduce $F(y)$ for an arbitrary $y\in TM$.  Note that every 
   nonzero $y$ is expressible as a positive multiple 
   $c$ of some $v$ with $F(v) = 1$. For $c > 0$, traversing 
   $c v$ under the windy conditions should take $c$ seconds. 
   Consequently, $F$ is positively homogeneous. Using this homogeneity and the
   formula derived for $F(v)$, we find that: 
   \begin{displaymath} 
    F(y)
    = 
      \frac{ \sqrt{ \, [h(W,y)]^2 + |y|^2 \lambda \, } }{\lambda} 
    - \frac{ \ h(W,y) \ }{\lambda} \, .
   \end{displaymath} 
   It is now manifest that $F(-y) \not= F(y)$. By hypothesis, $|W| < 1$, hence 
   $\lambda > 0$. We see from the formula for $F(y)$ that it is positive 
   whenever $y \not= 0$. Also, $F(0) = 0$ as expected.    

                 %...................................

  \subsubsection{New Riemannian metric and 1-form} 

   Our formula for $F$ has two parts. 
   \begin{itemize} 
    \item The first term is the norm of $y$ with respect to a new 
          Riemannian metric  
          \begin{displaymath} 
           a_{ij} =  \frac{h_{ij}}{\lambda}  
                   + \frac{ W_i}{\lambda} \frac{ W_j}{\lambda} \, ,  
          \end{displaymath} 
          where $W_i :=  h_{ij} \, W^j$ and $\lambda = 1 - W^i W_i$.  
    \smallskip 
    \item The second term is the value on $y$ of a differential 1-form 
          \begin{displaymath} 
           b_i = \frac{ - W_i }{\lambda} \, .  
          \end{displaymath} 
   \end{itemize} 

   Under the influence of $W$, the most efficient navigational paths are no 
   longer the geodesics of the Riemannian metric $h$; instead, they are the 
   geodesics of the Finsler metric $F$. For $\mathbb R^2$, this phenomenon is 
   treated by Carath\'eodory \cite{C99} as Zermelo's navigation 
   problem \cite{Z31}. Shen \cite{S02b} showed that the same phenomenon holds 
   for arbitrary Riemannian backgrounds in all dimensions.  

                %------------------------------------

 \subsection{Ubiquitous class of Finsler metrics} 

  The Finsler metric $F$ derived from the perturbation has the simple form 
  $F := \alpha + \beta$, where 
  \begin{displaymath}
   \alpha (x, y) := \sqrt{ \, a_{ij}(x) \, y^i y^j \, } \ , 
   \quad 
   \beta (x, y)  := b_i(x) \, y^i .   
  \end{displaymath}
  This is the defining feature of Randers metrics, which were introduced 
  by Randers in 1941 \cite{Ra41} in the context of general relativity, 
  and later named by Ingarden \cite{I57}. 

  The function $F$ is positive on the manifold $TM \setminus 0$, whose points 
  are of the form $(x,y)$, with $0 \not= y \in T_xM$. Over each point $(x,y)$ 
  of $TM \setminus 0$ (treated as a parameter space), we designate the vector 
  space $T_xM$ as a fiber, and name the resulting vector bundle $\pi^*TM$. 
  There is a canonical symmetric bilinear form $g_{ij} \, dx^i \otimes dx^j$ 
  on the fibers of $\pi^*TM$, with  
  \begin{displaymath} 
   g_{ij} := \half \left( F^2 \right)_{y^i y^j} \, . 
  \end{displaymath} 
  The subscripts $y^i, y^j$ signify partial differentiation, and the matrix 
  $(g_{ij})$ is known as the fundamental tensor. A Finsler metric $F$ 
  is said to be {\it strongly convex} if the said bilinear form is positive 
  definite, in which case it defines an inner product on each fiber 
  of $\pi^*TM$.  
  
  For a Randers metric to be strongly convex, it is necessary and sufficient 
  to have    
  \begin{displaymath} 
   \| b \| := \sqrt{ \, b_i \, b^i \, } < 1 \, , 
   \ \ \textrm{where} \ \   
   b^i := a^{ij} \, b_j \, .
  \end{displaymath}  
  See \cite{BCS00} or \cite{AIM93} for the proof of this fact. In our case, 
  using $a^{ij} = \lambda (h^{ij} - W^i W^j)$ and $b^i = - \lambda W^i$, we 
  find that  
  \begin{displaymath} 
   \| b \|^2 := a^{ij} \, b_i b_j = h_{ij} \, W^i W^j =: |W|^2 \, , 
  \end{displaymath} 
  which is less than $1$ by hypothesis. Therefore the described perturbation 
  of Riemannian metrics $h$ by vector fields $W$ with $|W| < 1$ always 
  generates strongly convex Randers metrics. 

                %------------------------------------
 
 \subsection{An inverse problem}   

  A question naturally arises: can every strongly convex Randers metric be 
  realized through the perturbation of some Riemannian metric $h$ by some 
  vector field $W$ satisfying $|W| < 1$ ? 

  Happily, the answer to this question is yes. Indeed, let us be given an 
  arbitrary Randers metric $F$ with data $a$ and $b$, 
  respectively a Riemannian metric and a differential 1-form, such that 
  $\| b \|^2 := a^{ij} \, b_i b_j < 1$.  Set $b^i := a^{ij} \, b_j$, and 
  $\varepsilon := 1 - \| b \|^2$. Construct $h$ and $W$ as follows: 
  \begin{displaymath} 
   h_{ij} := \varepsilon \, ( a_{ij} - b_i b_j ) \, , 
   \quad \quad    
   W^i := - b^i / \varepsilon \, .  
  \end{displaymath} 
 
  Note that $F$ is Riemannian if and only if $W=0$, in which case 
  $h = a$. Also, we have $W_i := h_{ij} \, W^j = - \varepsilon \, b_i$. 
  Using this, it can be directly checked that perturbing the above $h$ by the 
  stipulated $W$ gives back the Randers metric we started with. Furthermore,  
  \begin{displaymath} 
   |W|^2 := h_{ij} \, W^i W^j = a^{ij} \, b_i b_j =: \| b \|^2 < 1 \, . 
  \end{displaymath}   
  Let us summarize: 
   \begin{proposition}
    A strongly convex Finsler metric $F$ is of Randers type if and only if 
    it solves the Zermelo navigation problem on some Riemannian manifold 
    $(M, h)$, under the influence of a wind $W$ with $h(W,W) < 1$. 
    Also, $F$ is Riemannian if and only if $W=0$. 
   \end{proposition}

  Incidentally, the inverse of $h_{ij}$ is 
  $h^{ij} = \varepsilon^{-1} \, a^{ij} + \varepsilon^{-2} \, b^i b^j$.  This 
  $h^{ij}$, together with $W^i$, defines a Cartan metric $F^*$ of Randers 
  type on the cotangent bundle $T^*M$.  A comparison with \cite{HS96} shows 
  that $F^*$ is the Legendre dual of the Finsler-Randers metric $F$ on $TM$. 
  It is remarkable that the Zermelo navigation data of any strongly convex 
  Randers metric $F$ is so simply related to its Legendre dual.  
  See also \cite{Zi82} and \cite{S02b}.   

                 %-------------------------------------

 \subsection{Remark about isometries} 

  Two Finsler spaces $(M_1, F_1)$ and $(M_2, F_2 )$ are said to be 
  isometric if there exists a diffeomorphism $\phi : M_1 \to M_2$ which, when 
  lifted to a map between $TM_1$ and $TM_2$, satisfies $\phi^* F_2 = F_1$.

  Now consider two strongly convex Randers metrics $F_1$ and $F_2$, 
  where $F_i$ has Riemannian data $(a_i, b_i)$.  By the above proposition, 
  they arise as solutions to Zermelo's navigation problem with $(h_1, W_1)$ 
  and $(h_2,W_2)$, respectively.  A moment's thought gives the lemma below. 
  \begin{lemma}
   Let $\phi : M_1 \to M_2$ be a diffeomorphism. 
   The following three statements are equivalent: 
   \begin{itemize} 
    \item $\phi$ lifts to an isometry between $F_1$ and $F_2$. 
    \item $\phi^* a_2 = a_1$ and $\phi^* b_2 = b_1$. 
    \item $\phi^* h_2 = h_1$ and $\phi_* W_1 = W_2$. 
   \end{itemize} 
  \end{lemma}  

%------------------------------------------------------------------------------

\section{Zermelo navigation on Riemannian space forms} 

 This section illustrates a variety of perturbations on 3-dimensional 
 Riemannian space forms. In each example, with the exception of the radial 
 perturbation on the Euclidean metric (\S 3.2.2), $W$ is an infinitesimal 
 isometry of $h$. It happens that all the resulting strongly convex Randers 
 metrics are of constant flag curvature. The concept of flag curvature is a 
 natural extension of Riemannian sectional curvatures to the Finslerian realm; 
 see \S 3.4 for a review. 

 Since all our examples are in three dimensions, we let $(x,y,z)$ denote 
 position coordinates, and expand arbitrary tangent vectors as 
 $u \partial_x + v \partial_y + w \partial_z$. We give expressions for the 
 norm $\alpha := \sqrt{a(y,y)}$ instead of $a_{ij}$ because the former are 
 more compact. The Riemannian metric $a$ (defined in \S 2.2) can be recovered 
 via $a_{ij} = (\half \alpha^2)_{y^i y^j}$. 
  
                %------------------------------------

 \subsection{Spheres} 

  \subsubsection{Rotational perturbation}

   Let $S^3$ denote the standard unit sphere in $\mathbb R^4$. Using its 
   tangent spaces at the east and west poles, we may parametrize the sphere 
   by 
   \begin{displaymath}
    (x, y, z) \longrightarrow \frac{1}{ \sqrt{1 + x^2 + y^2 + z^2} }
                              (s, x, y, z) \, ; 
   \end{displaymath}
   here, $s=\pm 1$, respectively, for the eastern and western hemispheres. 
   Note that the equator corresponds to asymptotic infinity on the above 
   tangent spaces. Fix any constant $0 < \tau < 1$ and perturb via the 
   infinitesimal rotation 
   \begin{displaymath}
    W = \tau \, (y, -x, 0) \, ,  
    \quad \hbox{ with } \quad
    |W| = \tau \, \sqrt{\frac{x^2 + y^2}{1+x^2+y^2+z^2}} < 1 \, .  
   \end{displaymath}
   The bound on $\tau$ is needed to maintain $|W| < 1$ globally on $S^3$. The 
   resulting Randers metric $F = \alpha + \beta$ has constant flag curvature 
   $K=1$.  Explicitly,   
   \begin{eqnarray*}
    \alpha & = & \sqrt{
                 \frac{ \,  \rho^2 (u^2+v^2)  
                          - (\rho + \tau^2 \varphi) (xu+yv)^2 
                          + \eta \left\{  (\rho - z^2) w^2 
                                         - 2 zw (xu+yv) \right\} \, }
                      {\rho \, \eta^2} 
                      } \, ,                                                \\ 
    \beta  & = & \frac{\tau \, (-yu + xv)}{\eta} \, ,  
   \end{eqnarray*}
   where $\varphi := 1+z^2$, $\rho := 1+x^2+y^2+z^2$, 
   and $\eta := 1 + (1-\tau^2)(x^2+y^2) + z^2$. 

                %....................................

  \subsubsection{Perturbing by a privileged Killing field}

   Again, start with the unit sphere $S^3$ in $\mathbb R^4$, parametrized 
   as above. For each constant $K > 1$, let $h$ be $\textstyle{\frac{1}{K}}$ 
   times the standard Riemannian metric induced on $S^3$. The re-scaled metric 
   has sectional curvature $K$. 

   Perturb $h$ by the Killing vector field  
   \begin{displaymath} 
    W = \sqrt{K-1} \, \Big( -s(1+x^2), z-sxy, -y-sxz \Big) \, , 
    \quad \hbox{ with } \quad
    |W| = \sqrt{\frac{K-1}{K}} \, . 
   \end{displaymath} 
   This $W$ is tangent to the $S^1$ fibers in the Hopf fibration of $S^3$. 
   The resulting Randers metric $F$ has constant flag curvature $K$.  
   Moreover, it is not projectively flat \cite{BS02}.  This is in stark 
   contrast with the Riemannian case because, according to Beltrami's theorem, 
   a Riemannian metric is locally projectively flat if and only if it is of 
   constant sectional curvature.  Explicitly, $F = \alpha + \beta$, where
   \begin{eqnarray*} 
    \alpha &=& \frac{ \sqrt{ \, K(su-zv+yw)^2+(zu+sv-xw)^2+(-yu+xv+sw)^2 \, } }
                    {1+x^2+y^2+z^2} \, , 	                           \\
    \beta  &=& \frac{ \sqrt{K-1} \, (su-zv+yw) }{1+x^2+y^2+z^2} \, . 
   \end{eqnarray*} 

                %------------------------------------

 \subsection{Euclidean space}
   
  \subsubsection{Rotational perturbation}
  
   The Riemannian metric to be perturbed is the flat metric on $\mathbb R^3$. 
   The perturbing vector field is the infinitesimal rotation 
   $W := y \, \partial_x - x \, \partial_y + 0 \, \partial_z$. 
   The resulting Randers metric \cite{S02b} $F = \alpha + \beta$ solves the 
   least time problem for fish that are {\it surface}-feeding in a cylindrical 
   tank with a rotating current. $F$ is defined on the open cylinder 
   $x^2 + y^2 < 1$ in $R^3$, and has constant flag curvature $K=0$. 
   Explicitly,   
   \begin{eqnarray*}  
    \alpha &=& \frac{ \sqrt{ \, (-yu+xv)^2 + (u^2+v^2+w^2)(1-x^2-y^2) \, } } 
                    {1-x^2-y^2} \, ,                                    \\ 
    \beta  &=& \frac{-yu+xv}{1-x^2-y^2} \, , 
    \quad \hbox{ with } \quad |W|^2 = x^2+y^2 \, .      
   \end{eqnarray*} 

                %....................................
   
  \subsubsection{Radial perturbation}
  
   Again, we perturb the Euclidean metric, but this time $M$ is the open ball 
   of radius $R$ in $\mathbb R^3$, centered at the origin. The perturbing 
   vector field is the radial 
   $W = \tau (x \partial_x + y \partial_y + z \partial_z )$, 
   where $\tau$ is a constant. Impose the constraint 
   $| \tau | \leqslant \frac{1}{R}$   
   to ensure that $|W| < 1$ on $M$. The resulting Randers metric 
   $F = \alpha + \beta$ is of constant flag curvature $K = -\fourth \tau^2$, 
   and is given by 
   \begin{eqnarray*} 
    \alpha &=& \frac{ \sqrt{ \, \tau^2 (xu+yv+zw)^2  
                            + (u^2+v^2+w^2) \{1 - \tau^2 (x^2+y^2+z^2)\} \,} } 
                    {1 - \tau^2(x^2+y^2+z^2)} \, , 	              	   \\
     \beta &=& \frac{-\tau (xu+yv+zw)}
                    {1 - \tau^2 (x^2+y^2+z^2)} \, , 
     \quad \hbox{ with } \quad |W| = \sqrt{\tau^2 (x^2+y^2+z^2)} \, .
   \end{eqnarray*} 
   When $R = 1$ and $\tau = -1$, the perturbation generates the Funk metric 
   \cite{F} on the unit ball in $\mathbb R^3$. See also \cite{O83, S01}. 
   The Funk metric is isometric to the so-called Finslerian Poincar\'e ball. 
   A 2-dimensional version of the latter is analyzed in \cite{BCS00}. 

                %...................................

  \subsubsection{Perturbing by a translation} 

   As above, $h$ is the Euclidean metric $\delta_{ij}$. Choose any three 
   constants $p$, $q$, $r$ which satisfy $p^2+q^2+r^2 < 1$. We perturb $h$ by 
   the vector field 
   \begin{displaymath} 
    W = (p, q, r) \, , 
    \quad \hbox{ with } \quad |W| = \sqrt{p^2+q^2+r^2} \, .
   \end{displaymath}  
   The resulting Randers metric $F = \alpha + \beta$ has the form 
   \begin{eqnarray*} 
    \alpha &=& \frac{ \sqrt{ \, (pu+qv+rw)^2 
                              + (u^2+v^2+w^2) \{ 1 - (p^2+q^2+r^2) \} \,} }    
                    { 1 - (p^2+q^2+r^2) } \, ,                             \\  
    \beta  &=& \frac{ - (pu+qv+rw) } 
                    {1 - (p^2+q^2+r^2)} \, . 
   \end{eqnarray*} 
   This $F$ has constant flag curvature $K=0$, and is a (locally) Minkowski 
   metric.  

                %------------------------------------

 \subsection{Hyperbolic space} 
   
  \subsubsection{Rotational perturbation}

   Consider the Klein metric 
   \begin{displaymath}
    h_{ij} = \frac{(1-x^2-y^2-z^2)\delta_{ij} + x_i x_j}
                  {(1-x^2-y^2-z^2)^2}
   \end{displaymath}
   on the unit ball 
   $\mathbb{B}^3 := \{ (x,y,z) \in \mathbb{R}^3 : x^2+y^2+z^2 < 1\}$.
   Here $x_i := \delta_{is} x^s$.
   We perturb by the infinitesimal rotation
   \begin{displaymath}
    W = ( y , -x , 0)\, , 
    \quad \hbox{ with } \quad 
    |W| = \sqrt{\frac{x^2+y^2}{1-x^2-y^2-z^2}} \, .
   \end{displaymath}
   In order that $|W| < 1$, we restrict to the domain 
   $\{ 2x^2 + 2y^2 + z^2 < 1 \}$.  Define $\varphi = 1-2x^2-2y^2-z^2$.  
   Perturbing $h$ by $W$ produces a Randers metric $F=\alpha + \beta$, with  
   \begin{eqnarray*}
    \alpha & \! = \! & 
     \sqrt{\frac{\varphi \left[
                         (1-z^2)(u^2+v^2) + (1-x^2-y^2)w^2 + 2zw(xu+yv)
                        \right]
                 + (x^2+y^2)(yu-xv)^2}
                {(1-x^2-y^2-z^2) \varphi^2}}	\\
    \beta & \! = \! & \frac{-yu + xv}{\varphi} \, . 
   \end{eqnarray*}
   It is of constant flag curvature $K=-1$.

                 %.......................................

  \subsubsection{Perturbing by a type $(S) (\S 7.4.2)$ Killing field}

   In this section we will consider a more complicated perturbation 
   of the Klein metric.  The perturbing vector field is 
   \begin{displaymath}
    W = \tau \, ( 1-x^2 , z-xy , -y-xz ) \, .
   \end{displaymath}
   In order to effect $|W| < 1$, we restrict to the domain
   $\{ (1-\tau^2)x^2 + (1+\tau^2)(y^2+z^2) < 1-\tau^2 \}$.
   The resulting Randers metric $F = \alpha + \beta$ is also of constant flag 
   curvature $K = -1$. The Riemannian portion $\alpha$ of $F$ is substantial.
   Before writing it down we introduce some abbreviations: 
   $\psi_1 = 1-x^2-y^2-z^2$, 
   $\psi_2 = 1-x^2+y^2+z^2$, 
   $\psi_3 = \psi_1^2 + \tau^2 (-1-x^2+y^2+z^2)(y^2+z^2)$, 
   $\psi_4 = \psi_1^2 + \tau^2 (-1+x^2-y^2+z^2)(1-x^2)$, 	
   $\psi_5 = \psi_1^2 + \tau^2 (-1+x^2+y^2-z^2)(1-x^2)$, 
   $\psi_6 = 2 \tau^2 (z \psi_1 - xy \psi_2)$, 
   $\psi_7 = -2 \tau^2 (1-x^2) yz$, 
   $\psi_8 = - \tau^2 (y \psi_1 + xz \psi_2)$, 
   $\psi_9 = (1-\tau^2)(1-x^2) - (1+\tau^2)(y^2+z^2)$. 
   Then $\alpha$ is given implicitly by 
   \begin{displaymath}
    \psi_1 \psi_9^2 \, \alpha^2 
    =  
    \psi_1 (xu+yv+zw)^2       
    + \psi_3 u^2 + \psi_4 v^2 + \psi_5 w^2 
    + \psi_6 uv + \psi_7 vw + \psi_8 wu .
   \end{displaymath}
   The linear term is $\beta = \{ \tau (-u-zv+yw) \} / \psi_9$.

                   %....................................

  \subsubsection{Perturbing by a type $(T) (\S 7.4.3)$ Killing field}

   In our final perturbation of the Klein metric we consider the infinitesimal
   isometry 
   \begin{displaymath}
     W = ( 1 - y - x^2 , x - xy , -xz ) \, .
   \end{displaymath}
   The condition $|W|<1$ holds if we restrict to the domain 
   $(1-y)^2 < 1 - x^2 - y^2 - z^2$.  On this domain the perturbation produces 
   a strongly convex Randers metric $F = \alpha + \beta$ of constant flag 
   curvature $K = -1$.  For convenience, define 
   $\varphi_1 = 1 - x^2 - y^2 - z^2$, 
   $\varphi_2 = (1 + y)^2$, 
   $\varphi_3 = (1-y^2-z^2) \varphi_1 - x^2 \varphi_2$, 
   $\varphi_4 = (1-z^2) \varphi_1 - (1-x^2-z^2) \varphi_2$, 
   $\varphi_5 = (1-x^2-y^2) \varphi_9$, 
   $\varphi_6 = -2x (\varphi_1 + y \varphi_2)$, 
   $\varphi_7 = 2yz \varphi_9$, 
   $\varphi_8 = 2xz \varphi_9$, 
   $\varphi_9 = \varphi_1 - \varphi_2$.       
   We have 
   \begin{displaymath}
    \varphi_1 \varphi_9^2 \, \alpha^2 
    =  
      \varphi_3 u^2 + \varphi_4 v^2 + \varphi_5 w^2  
    + \varphi_6 uv + \varphi_7 vw + \varphi_8 wu \, ,
   \end{displaymath}
   and $\beta = \{ (xv - (y+1) u \} / \varphi_9$.

                %------------------------------------

 \subsection{Finsler metrics of constant flag curvature} 

  Given any Finsler metric $F$, the Chern connection on the pulled-back 
  tangent bundle $\pi^*TM$ gives rise to two curvature tensors, one of which, 
  $R_j{}^i{}_{kl}$, is analogous to the curvature tensor in Riemannian 
  geometry.  Indices on $R$ are raised and lowered by the fundamental tensor 
  $g_{ij}$ and its inverse $g^{ij}$. 

  At any point $x$ on $M$, a flag consists of a flagpole 
  $0 \not= y \in T_xM$ and a transverse edge $V \in T_xM$. The corresponding 
  {\it flag curvature} is defined as 
  \begin{displaymath} 
   K(x,y,V) 
   := 
   \frac{V^i \, ( y^j \, R_{jikl} \, y^l ) \, V^k} 
        {\ g(y,y) \, g(V,V) - [g(y,V)]^2 \ } \, .  
  \end{displaymath} 
  In the generic Finslerian setting, both the Chern $hh$-curvature $R$ 
  and the inner product $g$ (given by the fundamental tensor $g_{ij}$) 
  depend on the flagpole $y$. This dependence is absent whenever we specialize 
  to the Riemannian realm, in which case the flag curvature becomes the 
  familiar sectional curvature. For details and conventions, see \cite{BCS00}. 
  A Finsler metric is said to have constant flag curvature $K$ if $K(x,y,V)$ 
  has the constant value $K$ for all locations $x \in M$, flagpoles $y$, and 
  transverse edges $V$. 

  We note an interesting phenomenon shared by all our examples. In each case, 
  the constant flag curvature of the resulting Randers metric $F$ does not 
  exceed the constant sectional curvature of the original Riemannian 
  metric $h$. 

%-----------------------------------------------------------------------------

\section{Navigation description of Randers metrics of constant flag curvature} 
 
                %------------------------------------
  
 \subsection{Characterization}  

  Let $F = \alpha + \beta$, with $\alpha^2 := a_{ij} \, y^i y^j$ and 
  $\beta := b_i \, y^i$, be a Randers metric.  Using $a^{ij}$ to raise the 
  index on the components $b_j$ of the 1-form $b$, we get a vector field 
  $b^\sharp = b^i \partial_{x^i}$.  Let us introduce the abbreviations  
  \begin{displaymath}
   \txtcurl_{ij} := \partial_{x^j} b_i - \partial_{x^i} b_j  
   \quad \quad  \textrm{and} \quad \quad 
   \theta_j := b^i \, \txtcurl_{ij} \, . 
  \end{displaymath} 
  Note that $\txtcurl$ is the 2-form $- d b$, and interior multiplication of 
  $\txtcurl$ by the vector field $b^\sharp$ gives the 1-form $\theta$. 

  Define the geometric quantity 
  \begin{displaymath} 
   \sigma := \frac{2 \, \textrm{div} \, b^\sharp}{ \ n - \| b \|^2 \ }  \, . 
  \end{displaymath} 
  A theorem in \cite{BR03} states that the Randers metric $F$ has constant 
  flag curvature $K$ if and only if $\sigma$ is constant,  
  \begin{displaymath}
   \mathcal L_{b^\sharp} a 
   = \sigma (a - b \otimes b) - (b \otimes \theta + \theta \otimes b) 
  \end{displaymath}
  (where 
   $\mathcal L_{b^\sharp} a
    =   b^k \, \partial_{x^k} a_{ij} 
      + a_{kj} \, \partial_{x^i} b^k 
      + a_{ik} \, \partial_{x^j} b^k$ 
  is a Lie derivative), and the Riemann tensor of $a$ has the form 
  \begin{eqnarray*}  
   \aR_{hijk}                                                       
   &=&   \  \xi \, ( a_{ij} \, a_{hk} - a_{ik} \, a_{hj} )       \\ 
   & & - \, \fourth \, a_{ij} \, \txtcurl^t{}_h \, \txtcurl_{tk} 
       +    \fourth \, a_{ik} \, \txtcurl^t{}_h \, \txtcurl_{tj} \\ 
   & & + \, \fourth \, a_{hj} \, \txtcurl^t{}_i \, \txtcurl_{tk} 
       -    \fourth \, a_{hk} \, \txtcurl^t{}_i \, \txtcurl_{tj} \\ 
   & & - \, \fourth \, \txtcurl_{ij} \, \txtcurl_{hk} 
       +    \fourth \, \txtcurl_{ik} \, \txtcurl_{hj}                     
       +    \half \, \txtcurl_{hi} \, \txtcurl_{jk} \, ,         \\
   \textrm{with} \quad 
   \xi \!\!\!\! &:=& \!\!\!\! 
                (K - \threesixteenths \sigma^2) 
              + (K + \onesixteenth \sigma^2) \, \| b \|^2 
              - \fourth \, \theta^i \theta_i \, .  
  \end{eqnarray*}    
  In these equations, all tensor indices are raised and lowered by $a$. 
  For later purposes, let us refer to the above as the Basic equation and 
  the Curvature equation, respectively.

  The Basic equation {\it alone} is equivalent to the statement that the 
  S-curvature (divided by $F$) has the constant value $\fourth \sigma (n+1)$; 
  see \cite{CS03}.  While the Basic equation only makes sense for Randers 
  metrics, its characterization in terms of the S-curvature gives a 
  well-defined criterion which can be imposed on Finsler metrics in general. 

                %-------------------------------------

 \subsection{Matsumoto identity}   

  In the original statement of the characterization above, there is a third 
  equation, named CC(23), that $a$ and $b$ must satisfy.  As such, the said 
  theorem is equivalent in content to one in \cite{MS02}.  Recent work shows 
  that the CC(23) equation is derivable from the Basic and Curvature 
  equations with $\sigma$ constant.  Hence it is omitted here.  See 
  \cite{BR04} for more discussions.  

  The omitted CC(23) equation is geometrically significant because, in the 
  presence of a preliminary form of the Basic and Curvature equations, it is 
  equivalent to the constancy of $\sigma$ (or the S-curvature).  The CC(23) 
  equation is also useful.  For example, it leads to the 
  following Matsumoto identity, which describes the interplay between 
  $\sigma$ and $K$: 
  \begin{displaymath} 
   \sigma (K + \onesixteenth \sigma^2) = 0 
   \quad \textrm{for constant} \ K \ \textrm{and} \ n \geqslant 2 .  
  \end{displaymath}   
   
                %--------------------------------------

 \subsection{Navigation description}  
 
  According to Proposition 1, our strongly convex Randers metric $F$ can be 
  realized as the perturbation of a Riemannian metric $h$ by a vector field 
  $W$ which satisfies $h(W,W) < 1$. Using this fact and \S 2.1.3, the tensors 
  $a$ and $b$ that comprise $F$ are expressible as 
  \begin{displaymath} 
   a_{ij} =  \frac{h_{ij}}{\lambda}  
           + \frac{ W_i}{\lambda} \frac{ W_j}{\lambda} \, ,  
   \quad \quad   
   b_i = \frac{ - W_i }{\lambda} \, ,  
  \end{displaymath}  
  where $W_i :=  h_{ij} \, W^j$ and $\lambda := 1 - h(W,W) > 0$.  For $a^{ij}$ 
  and $b^i$, see \S 2.2. 

                %......................................

  \subsubsection{Navigation version of the Basic equation} 

   The Basic equation in the stated characterization involves $a$, $b$, 
   $\mathcal L_{b^\sharp} a$, and $\theta$.  Substituting the above formulae 
   for $a$, $b$ and computing the requisite partial derivatives in the 
   remaining two tensors, we obtain an equivalent $\mathcal L_W$ equation:   
   \begin{displaymath} 
    \mathcal L_W h = - \sigma \, h \, .  
   \end{displaymath} 
   The left-hand side can be rewritten in terms of the covariant derivative 
   operator ``${}_{:}$" associated to $h$, and the $\mathcal L_W$ equation 
   becomes 
   \begin{displaymath} 
    W_{i:j} + W_{j:i} = - \sigma \, h_{ij} \, . 
   \end{displaymath} 
   In this equation, 
   \begin{center} 
    ``$\sigma$ must vanish whenever $h$ is not flat." 
   \end{center}    
   Indeed, let $\varphi_t$ denote the time $t$ flow of the vector field $W$. 
   The $\mathcal L_W$ equation tells us that 
   $\varphi_t^* h = e^{- \sigma t} h$. Since $\varphi_t$ is a diffeomorphism, 
   $e^{- \sigma t} h$ and $h$ must be isometric; therefore they have the same 
   sectional curvatures. If $h$ is not flat, this condition on sectional 
   curvatures mandates that $e^{- \sigma t} = 1$, hence $\sigma = 0$.  The 
   above argument was pointed out to us by Bryant.   

                  %.....................................

  \subsubsection{Riemannian connections of $a$ and $h$} 

   To minimize some anticipated clutter, let us introduce the abbreviations  
   \begin{displaymath} 
    \mathcal C_{ij} := \partial_{x^j} W_i - \partial_{x^i} W_j 
                     = W_{i:j} - W_{j:i} \, , 
    \quad \quad 
    \mathcal T_j := W^i \, \mathcal C_{ij} \, , 
   \end{displaymath} 
   and agree to let the subscript $0$ denote contraction of any index with 
   $y^i$. Indices on $\mathcal C$, $\mathcal T$ are to be manipulated by the 
   Riemannian metric $h$ only. 

   Let $\agamma^i{}_{jk}$ and $\aG^i := \half \agamma^i{}_{00}$ be, 
   respectively, the Christoffel symbols and geodesic spray coefficients of 
   the Riemannian metric $a$.  Likewise, let $\hG^i := \half \hgamma^i{}_{00}$ 
   be the geodesic spray coefficients of $h$.  (The factor of $\half$ here is 
   absent in some references such as \cite{BCS00}.)  A straight-forward 
   computation, or an application of Rapcs\'ak's identity \cite{Rap61}, 
   together with the $\mathcal L_W$ equation, shows \cite{BR04} that  
   \begin{displaymath} 
    \aG^i 
    =   \hG^i 
      + \frac{y^i}{ \ 2 \lambda \ }(\mathcal T_0 - \sigma \, W_0)
      - \mathcal T^i \left(  \frac{h_{00}}{ \ 4 \lambda \ }  
                           + \frac{W_0 \, W_0}{ \ 2 \lambda^2 \ } \right) 
      + \frac{ \ \mathcal C^i{}_0 \, W_0 \ }{2 \lambda} \, . 
   \end{displaymath} 

                      %................................
 
  \subsubsection{Navigation version of the Curvature equation} 

   Abbreviate the above formula as $\aG^i = \hG^i + \zeta^i$.  We now use it 
   to relate the curvature tensor $\aR$ of $a$ to the curvature tensor $\hR$ 
   of $h$. To this end, consider the spray curvature \cite{B47b} tensors 
   $\aK^i{}_j = \aR_0{}^i{}_{j0}$ and $\hK^i{}_j = \hR_0{}^i{}_{j0}$.  The 
   Riemann tensor can be recovered from the spray curvature through  
   $\aR_{hijk} = \third \{ (\aK_{ij})_{y^k y^h} - (\aK_{ik})_{y^j y^h} \}$,  
   where the up index on $\aK$ has been lowered by $a$.  A similar formula 
   holds for $\hR_{hijk}$ and $\hK_{ij}$, with the index on $\hK$ lowered by 
   $h$.  The advantage of working with the spray curvature is that it has 
   less indices than the full Riemann tensor.  

   The Curvature equation of \S 4.1 can be recast into the form   
   \begin{eqnarray*} 
    \aK^i{}_j  
    &=  &\xi \, (\alpha^2 \, \delta^i{}_j - y^i \, \ay_j)                 \\ 
    &   &+ \, \fourth \, \txtcurl^s{}_0 \, 
              (  \txtcurl_s{}^i \, \ay_j 
               + y^i \, \txtcurl_{sj} 
               - \txtcurl_{s0} \, \delta^i{}_j )                          \\   
    &   &- \, \fourth \, \alpha^2 \, \txtcurl^{si} \, \txtcurl_{sj} 
            - \threefourths \, \txtcurl^i{}_0 \, \txtcurl_{j0} \, , 
   \end{eqnarray*} 
   where $\xi$ is as defined in \S 4.1 and $\ay_j := a_{jk} y^k$. 
   Into (the left-hand side of) this we substitute one version of the split 
   covariantized Berwald formula (see \cite{BR04, S01} for expositions and 
   references therein), which says that 
   \begin{displaymath} 
    \aK^i{}_j 
    =   \hK^i{}_j 
      + (2 \, \zeta^i)_{:j} - (\zeta^i)_{y^s} (\zeta^s)_{y^j}
      - y^s (\zeta^i{}_{:s})_{y^j} + 2 \, \zeta^s (\zeta^i)_{y^s y^j} \, . 
   \end{displaymath} 
   Here, the subscripts ``${}_{y^k}$" mean $\partial_{y^k}$. 
   This is followed by a tedious calculation, in which all quantities are 
   rewritten in terms of the navigation variables $h$, $W$, and the 
   $\mathcal L_W$ equation is used prodigiously. A formula for 
   $\hK^i{}_j$ then results, from which we compute the Riemann tensor 
   $\hR_{hijk}$. 
    
   The outcome of that calculation is remarkable. It says that given the  
   $\mathcal L_W$ equation, the said Curvature equation is equivalent to the 
   statement that $h$ is a Riemannian space form of constant sectional 
   curvature $K + \onesixteenth \sigma^2$. Namely, 
   \begin{displaymath} 
    \hR_{hijk} = ( K + \onesixteenth \sigma^2 ) 
                 (h_{ij} \, h_{hk} - h_{ik} \, h_{hj}) \, . 
   \end{displaymath}  

                %----------------------------------

 \subsection{Summary} 

  \begin{theorem} 
   A strongly convex Randers metric $F$ has constant flag curvature $K$ if 
   and only if:  
   \begin{itemize} 
    \item $F$ solves Zermelo's navigation problem on a Riemannian space form 
          $(M,h)$ of sectional curvature $K + \onesixteenth \sigma^2$ for some 
          constant $\sigma$, under the influence of a vector field 
          (``wind") $W$.  
    \item The wind $W$ satisfies $h(W,W) < 1$, and is coupled to $h$ and 
          $\sigma$ in such a way that 
          $\mathcal L_W h = - \sigma \, h$, where $\mathcal L$ denotes 
          Lie differentiation. 
   \end{itemize} 
   For non-flat $h$, $\sigma$ must vanish, in which case $W$ must be a 
   Killing vector field of $h$. 
  \end{theorem} 

  The last statement has already been observed in \S 4.3.1.  Alternatively, 
  since the sectional curvature of $h$ is $K + \onesixteenth \sigma^2$, that  
  statement also follows from Matsumoto's identity (\S 4.2). 

  Note that $K$, the flag curvature of $F$, is bounded above by the 
  sectional curvature $K + \onesixteenth \sigma^2$ of $h$. This explains the 
  phenomenon we noted at the end of \S 3.4.  Since $\sigma$ must vanish 
  whenever $K + \onesixteenth \sigma^2 \not= 0$, we have the following 
  trichotomy. 
  \begin{itemize} 
   \item[$(+)$] 
        For $K > 0$: The quantity $K + \onesixteenth \sigma^2$ is positive, 
        hence $\sigma = 0$. Consequently the sectional curvature of $h$ must 
        equal $K$, the flag curvature of $F$. 
   \item[$(0)$] 
        For $K = 0$: The sectional curvature of $h$ reduces to 
        $\onesixteenth \sigma^2$. If $\sigma$ were nonzero, $h$ would have to 
        be flat according to the last part of Theorem 3; but that would be 
        incompatible with having sectional curvature $\onesixteenth \sigma^2$. 
        So $\sigma$ must vanish, whence $h$ is flat. 
   \item[$(-)$] 
        For $K < 0$: There are two viable scenarios. The first is 
        $\sigma = \pm 4 \sqrt{|K|}$, in which case $h$ is flat. For the second 
        scenario, $K + \onesixteenth \sigma^2 \not= 0$; hence $\sigma = 0$ 
        and $h$ must have negative sectional curvature $K$. 
  \end{itemize} 
  
%------------------------------------------------------------------------------

\section{Complete list of allowable vector fields} 

  Our goal here is towards a classification of Randers metrics of constant 
  flag curvature. By the navigation description, these metrics arise as 
  perturbations of Riemannian space forms $h$ by vector fields $W$ satisfying 
  $W_{i:j} + W_{j:i} = -\sigma \, h_{ij}$. For each of the three standard 
  models (Euclidean, spherical and hyperbolic) of Riemannian space forms we 
  derive a formula for $W$. 

                %------------------------------------  

 \subsection{Setting some notation with a basic lemma} 

   \begin{lemma}\label{lem4}
    Let $P_i = P_i(x)$ be solutions of the following system 
    \begin{displaymath} 
     \frac{\pa P_i}{\pa x^j} + \frac{\pa P_j}{\pa x^i} = 0 . 
    \end{displaymath} 
    Then 
    \begin{displaymath} 
     P_i = Q_{ij} \, x^j + C_i \, , 
    \end{displaymath} 
    where $(C_i)$ is an arbitrary constant row vector and $Q = (Q_{ij})$ is 
    an arbitrary constant skew-symmetric matrix ($Q_{ji} = - Q_{ij}$).  
   \end{lemma}
  
   {\it Proof} : Using the defining differential equation three times, we have 
   \begin{displaymath}
       \frac{\partial^2 P_i}{\partial x^k \partial x^j} 
    =  
     - \frac{\partial^2 P_j}{\partial x^k \partial x^i}			
    =  
       \frac{\partial^2 P_k}{\partial x^i \partial x^j}			
    =  
     - \frac{\partial^2 P_i}{\partial x^j \partial x^k} \, .
   \end{displaymath}
   This shows that all second order partial derivatives of $P_i$ must vanish. 
   Hence $P_i$ must be linear; that is, it has the form 
   $P_i = Q_{ij} x^j + C_i$, with constants $Q_{ij}$ and $C_i$. Inserting this 
   expression into the defining PDE shows that $Q_{ij} + Q_{ji} = 0$. \qed

   \begin{quote} 
    For the rest of the paper: ``$\cdot$'' refers to the standard dot 
    product on $\mathbb{R}^n$; indices on $Q$ and $C$ are raised and lowered 
    by the Kronecker delta $\delta_{ij}$; and $Qx+C$ means 
    $(Q^i{}_j \, x^j + C^i)$.  We regard $(C_i)$ as a row vector and $(C^i)$ 
    as a column vector. 
   \end{quote} 

                  %...................................

 \subsection{The Euclidean case} 

   The first Riemannian space form we consider is the flat Euclidean metric.  
   The admissible perturbing vector fields $W$ are described in the 
   following proposition. 

   \begin{proposition} 
    Let $F = \alpha + \beta$ be a strongly convex Randers metric which results 
    from perturbing the flat metric $h_{ij} = \delta_{ij}$ on $\mathbb R^n$ by 
    a vector field $W = (W^i)$. Then $F$ is of constant flag curvature $K$ if 
    and only if $W$ has the form
    \begin{displaymath} 
     W^i(x) = - \half \sigma \, x^i + Q^i{}_j \, x^j + C^i \, ,     
    \end{displaymath} 
    where $(Q^i{}_j)$ is a constant skew-symmetric matrix, $(C^i)$ is a 
    constant column vector, $\sigma$ is a constant such that 
    $\sigma^2 = - 16 K$, and 
    \begin{displaymath} 
       (Qx+C) \cdot (Qx+C) 
     + \sigma x \cdot \left ( \fourth \sigma x - C \right ) 
     \ < \ 1 \, . 
    \end{displaymath}   
   \end{proposition}

   \noindent{Remark:} Note that by virtue of $\sigma^2 = - 16 K$, we see that 
   $K$ must be $\leqslant 0$. 

   {\it Proof} : 
   Being flat, $h$ satisfies the space form criterion of the navigation  
   description, with $K + \onesixteenth \sigma^2 = 0$. The rest of the proof 
   studies the second criterion, which is the equation 
   $\mathcal L_W h = - \sigma \, h$.  

   \begin{itemize} 
    \item[($\Leftarrow$)] 
         Suppose $W$, with its index lowered by $h_{ij} = \delta_{ij}$, is of 
         the form
         \begin{displaymath} 
          W_i = - \half \sigma \, \delta_{ij} \, x^j 
                + Q_{ij} \, x^j + C_i \, . 
         \end{displaymath} 
         Keeping in mind that the covariant derivative ``${}_:$" associated 
         with the Euclidean $h$ is simply partial differentiation, together 
         with the skew-symmetry of $Q$, we immediately obtain
         \begin{displaymath} 
          ({}^*) \quad \quad \quad \quad 
          W_{i:j} + W_{j:i} = - \sigma \, \delta_{ij} \, . 
         \end{displaymath} 
         Thus the $\mathcal L_W$ equation in the navigation description is 
         satisfied, and $F$ has constant flag curvature $K$. 
    \smallskip 
    \item[($\Rightarrow$)] 
         Conversely, suppose $F$ has constant flag curvature $K$.
         By the navigation description, $W$ must be a solution of (*).
         Note that 
         \begin{displaymath} 
          W_i = - \half \sigma \, \delta_{ij} \, x^j 
         \end{displaymath} 
         is a particular solution. Adding to it the solutions of the 
         homogeneous system
         $\frac{\pa P_i}{\pa x^j} + \frac{\pa P_j}{\pa x^j} = 0$ 
         gives the general solution. According to Lemma \ref{lem4}, the 
         latter have the form $P_i = Q_{ij} \, x^j + C_i$, where each $C_i$ 
         is constant and $(Q_{ij})$ is a constant skew-symmetric matrix. 
         Using $h^{ij} = \delta^{ij}$, we raise the index on $W_i$ to effect 
         the $W^i$ as claimed. 
   \end{itemize}  
   The inequality satisfied by $Q$, $C$, and $\sigma$ comes from the 
   requirement $|W| < 1$. \qed

                     %..................................

 \subsection{The spherical and hyperbolic cases} 

   We now perturb standard models of Riemannian metrics with constant 
   sectional curvature $\kappa \not= 0$. The list of allowable $W$ is given 
   in the following proposition. 
 
   \begin{proposition} 
    Let $F= \alpha+\beta$ be a strongly convex Randers metric which results 
    from perturbing the standard, complete, simply connected, $n$-dimensional 
    Riemannian space $(M,h)$ of constant sectional curvature $\kappa \not= 0$ 
    by a vector field $W$. Then $F$ is of constant flag curvature $K$ if and 
    only if $K = \kappa$ and $W$ is Killing, with the following description in 
    terms of a constant vector $(C^i)$ and a constant skew-symmetric matrix 
    $(Q^i{}_j)$.
    \begin{enumerate} 
     \item[(a)] $K = \kappa > 0$. Employ a projective coordinate system on 
                the unit $n$-sphere, one which comes from parametrizing each 
                hemisphere using the tangent space at the pole.  Multiply the 
                standard Riemannian metric by $\frac{1}{K}$ to effect constant 
                sectional curvature $K$.  The $h$-norm of any tangent vector 
                $y$ is given by
                \begin{displaymath} 
                 \quad 
                 |y| := \sqrt{h(y,y)}  
                 = \frac{1}{ \ \sqrt{K} \ }
                 \frac{ \sqrt{  (y \cdot y)(1 + x \cdot x) - (x \cdot y)^2 } } 
                      {1+ x \cdot x} \, , 
                 \quad  y \in T_x \mathbb R^n \simeq \mathbb R^n . 
                \end{displaymath}              
                With respect to this coordinate system, 
                \begin{displaymath}  
                 W^i(x) 
                 = 
                 Q^i{}_j \, x^j + C^i + (x \cdot C) \, x^i . 
                \end{displaymath} 
     \item[(b)] $K = \kappa < 0$. Let $h$ be the Klein model of constant  
                sectional curvature $K$ on the unit ball ${\rm B}^n$, with the 
                Cartesian coordinates of $\mathbb R^n$. 
                The $h$-norm of any tangent vector $y$ is given by
                \begin{displaymath}              
                 \quad \quad 
                 |y| := \sqrt{h(y,y)}   
                 = \frac{1}{ \ \sqrt{|K|} \ }
                 \frac{ \sqrt{  (y \cdot y)(1 - x \cdot x) + (x \cdot y)^2 } } 
                      {1- x \cdot x} \, , 
                 \quad  y \in T_x \mathbb R^n \simeq \mathbb R^n . 
                \end{displaymath} 
                With respect to this coordinate system,              
                \begin{displaymath} 
                 W^i(x) 
                 = 
                 Q^i{}_j \, x^j + C^i - (x \cdot C) \, x^i .
                \end{displaymath} 
    \end{enumerate} 
    In each case, $W$ is subject to the constraint
    \begin{displaymath} 
     \frac{1}{ \, 1 + \psi (x \cdot x) \, } 
     \{ (Qx+C) \cdot (Qx+C) + \psi (x \cdot C)^2 \} 
     \ < \ |K| \, , 
     \quad \textrm{where} \ 
     \psi := \frac{K}{|K|} \, . 
    \end{displaymath}
   \end{proposition}
  
   {\it Proof} : Our Riemannian metric $h$ has constant sectional curvature 
   $\kappa \not=0$. Therefore it satisfies the space form criterion of 
   the navigation description, with $K + \frac{\sigma^2}{16} = \kappa$. 
   In particular, $K + \frac{\sigma^2}{16} \not= 0$. The Matsumoto identity 
   then implies that $\sigma$ must vanish. Consequently, $K = \kappa$. 
   
   According to our navigation description, perturbing the above $h$ by 
   a vector field $W$ (with $|W|<1$) generates a Randers metric of constant 
   flag curvature $K$ if and only if the equation 
   $\mathcal L_W h = - \sigma \, h$ is satisfied. Since $\sigma = 0$ here, 
   that equation reduces to the statement that $W$ is a Killing vector field 
   of $h$. The proof of this proposition therefore concerns the classification 
   of solutions of the Killing field equation: 
   \begin{displaymath} 
    W_{i:j}+ W_{j:i} = 0 \, .
   \end{displaymath} 
   
   \begin{itemize} 
    \item To minimize notational clutter, let us introduce the abbreviations  
          \begin{displaymath} 
           x_i  := \delta_{ij} \, x^j \, , 
           \quad \quad                 
           \rho := 1 + \psi (x \cdot x) .
          \end{displaymath}  
          Then
          \begin{displaymath} 
           h_{ij} 
           = 
           \frac{1}{ \ |K| \ }
           \left(
           \frac{\delta_{ij}}{\rho} 
           - \psi \frac{ \ x_i x_j \ }{\rho^2 } 
           \right) \, ,  
           \quad \quad  
           h^{ij} 
           = 
           \rho \, |K| \, \{ \delta^{ij} + \psi \, x^i x^j \} .
          \end{displaymath} 
          The Christoffel symbols of $h$ are given by
          \begin{displaymath} 
           \hgamma^k{}_{ij} 
           = - \psi \frac{ \ x_i \, \delta^k{}_j 
                           + x_j \, \delta^k{}_i \ }
                         {\rho} . 
          \end{displaymath} 
          Hence 
          \begin{displaymath} 
           W_{i:j} 
           =   \frac{\pa W_i}{\pa x^j} 
             + \psi \frac{ \ x_i W_j + x_j W_i \ }{\rho} . 
          \end{displaymath} 
          The Killing field equation now reads  
          \begin{displaymath}  
             \frac{\pa W_i}{\pa x^j} + \frac{\pa W_j}{\pa x^i} 
           + \frac{ \ 2 \psi \ }{\rho} (x_i W_j + x_j W_i) = 0 \, . 
          \end{displaymath} 
    \item To solve it, let us replace the dependent variables $W_i$ by new 
          ones that are named $P_i$, as follows:    
          \begin{displaymath}   
           W_i = \frac{1}{\ \rho \, |K| \ } P_i \, .
          \end{displaymath} 
          (The division by $|K|$ effects a simplification later, when we 
          use $h^{ij}$ to raise the index on $W_i$.) Computations give: 
          \begin{eqnarray*} 
           \quad \quad 
           \frac{\pa W_i}{\pa x^j} + \frac{\pa W_j}{\pa x^i} 
           &=&  
           \frac{1}{ \ \rho \, |K| \ }
           \left ( \frac{\pa P_i}{\pa x^j} + \frac{\pa P_j}{\pa x^i} \right ) 
           - \frac{\ 2 \psi \ }{ \ \rho^2 \, |K| \ }
             (x_i P_j + x_j P_i) \, ,                                       \\ 
           \quad \quad 
           \frac{\ 2 \psi \ }{\rho} (x_i W_j + x_j W_i) 
           &=& 
           \frac{\ 2 \psi \ }{ \ \rho^2 \, |K| \ }
                (x_i P_j + x_j P_i) \, . 
          \end{eqnarray*} 
          This change of dependent variables transforms the above equation  
          into 
          \begin{displaymath} 
           \frac{\pa P_i}{\pa x^j} + \frac{\pa P_j}{\pa x^i} = 0 . 
          \end{displaymath} 
          By Lemma \ref{lem4}, the solutions $P_i$ have the form 
          \begin{displaymath}  
           P_i = Q_{ij} \, x^j + C_i \, , 
          \end{displaymath} 
          where $(Q_{ij})$ is a constant skew-symmetric matrix, and the $C_i$ 
          are constants. 
   \end{itemize} 
   Thus the covariant form (that is, with index down) of the Killing 
   field $W$ is   
   \begin{displaymath} 
    W_i = \frac{ \ Q_{ij} \, x^j + C_i \ }{\rho \, |K| } . 
   \end{displaymath} 
   To obtain the contravariant form (namely, with index up) of $W$, we raise 
   its index using 
   $h^{ij} = \rho \, |K| \, \{ \delta^{ij} + \psi \, x^i x^j \}$. 
   The result reads:  
   \begin{displaymath}  
    W^i 
    := 
    h^{ij} W_j 
    =  
    Q^i{}_j \, x^j + C^i + \psi (x \cdot C) \, x^i \, , 
   \end{displaymath} 
   where $Q^i{}_j := \delta^{is} Q_{sj}$ and $C^i := \delta^{is} C_s$. 
   
   Finally, the constraint on $Q$ and $C$ comes from the requirement 
   $|W| < 1$.     
   \qed

                %....................................

 \subsection{Remarks} 

  Note that: in the case of flat $h$, both $W_i$ and $W^i$ are polynomials of 
  degree 1 in the position variables $x$; for non-flat $h$, $W_i$ is a 
  rational function in $x$ of degree -1, while $W^i$ is a polynomial of 
  degree 2 in $x$ whenever $C \not= 0$. 
 
  We tabulate below the constant skew-symmetric matrix $Q$, the constant 
  vector $C$, and the value of the constant $\sigma$, for all the examples 
  of \S 3.  To reduce clutter, let $0_{3 \times 3}$ denote the 
  $3$-by-$3$ zero matrix, and 
  \begin{displaymath}
   J := \left(
   \begin{array}{rr}
   0  & 1 \\
   -1 & 0 \\ 
   \end{array}
   \right) \, .
  \end{displaymath}    

  \medskip 

  \begin{displaymath}
   \begin{array}{cccc}
   \hbox{Example} & (Q_{ij})              & (C_i)               & \sigma   \\
   \hline 
                  &                       &                     &          \\
   3.1.1          & \tau J \oplus 0       & (0,0,0)             & 0        \\ 
                  &                       &                     &          \\
   3.1.2          & 0\oplus\sqrt{K-1}\, J & (-s \sqrt{K-1},0,0) & 0        \\ 
                  &                       &                     &          \\
   3.2.1          & J \oplus 0            & (0,0,0)             & 0        \\ 
                  &                       &                     &          \\
   3.2.2          & 0_{3 \times 3}        & (0,0,0)             & -2 \tau  \\ 
                  &                       &                     &          \\
   3.2.3          & 0_{3 \times 3}        & (p,q,r)             & 0        \\ 
                  &                       &                     &          \\
   3.3.1          & J \oplus 0            & (0,0,0)             & 0        \\
                  &                       &                     &          \\ 
   3.3.2          & 0 \oplus \tau J       & (\tau,0,0)          & 0        \\
                  &                       &                     &          \\
   3.3.3          & J \oplus 0            & (1,0,0)             & 0 
   \end{array}
  \end{displaymath}

%------------------------------------------------------------------------------

\section{Classification of Randers metrics with constant flag curvature} 

 \subsection{The main theorem} 
 
  We now combine the navigation description (see \S 4.4) and the work of 
  \S 5 to classify Randers metrics of constant flag curvature. Before stating 
  the theorem, we recall that:
  \begin{itemize} 
   \item the skew-symmetric matrix $Q = (Q^i{}_j)$ and the vector $C = (C^i)$ 
         are constant; 
   \item $Qx$ denotes $(Q^i{}_j \, x^j)$, and $x := (x^i)$; 
   \item all indices on $Q$, $C$, $x$ are manipulated by the Kronecker 
         deltas $\delta_{ij}$ and $\delta^{ij}$;  
   \item ``$\cdot$'' is the standard Euclidean dot product. 
  \end{itemize} 
 
  \begin{theorem}[Classification]
   Let $F(x,y) =  \sqrt{a_{ij}(x) \, y^i \, y^j } + b_i(x) \, y^i$ 
   be a strongly convex Randers metric on a smooth manifold $M$ of dimension 
   $n \geqslant 2$.  Then $F$ is of constant flag curvature $K$ if and only if
   the following conditions are satisfied.
 
   \medskip 
   \noindent ${\bf (1)}$   
        The Riemannian metric $a$ and 1-form $b$ have the representation 
        \begin{displaymath}
         a_{ij} =  \frac{h_{ij}}{\lambda}  
                 + \frac{ W_i}{\lambda} \frac{ W_j}{\lambda} \, ,  
         \quad \quad 
         b_i = \frac{ - W_i }{\lambda} \, ,
        \end{displaymath}
        where $h$ is a Riemannian space form and $W = W^i \partial_{x^i}$ is an
        infinitesimal homothety (of $h$), both globally defined on $M$. Here, 
        $W_i := h_{ij} W^j$ and $\lambda := 1 - h(W,W) > 0$. 

   \medskip 
   \noindent ${\bf (2)}$     
        Up to local isometry, the Riemannian space form $h$ and the vector 
        field $W$ must belong to one of the following four families. 
        \begin{itemize}
         \item[$(+)$] 
              When {\bf $K > 0$:} 
              $h$ is $\frac{1}{K}$ times the standard metric on the unit 
              $n$-sphere, and $W = Qx + C + (x \cdot C) \, x$, with 
	      \begin{displaymath}
               \frac{1}{\ 1 + (x \cdot x) \ }
               \{ (Qx+C) \cdot (Qx+C) + (x \cdot C)^2 \} 
               \ < \ K \, .
              \end{displaymath} 
        \item[$(0)$] 
             When {\bf $K = 0$:}  
             $h$ is the Euclidean metric on $\mathbb R^n$ and $W = Qx+C$, with
             \begin{displaymath}  
              (Qx+C) \cdot (Qx+C) \ < \ 1 \, .   
             \end{displaymath} 
        \item[$(-)$] 
             When {\bf $K < 0$:}      
             \begin{itemize}
              \item[$(-)_e$] 
                    either $h$ is the Euclidean metric on $\mathbb R^n$, and 
                    $W = - \half \sigma x + Qx+C$ satisfies   
                    \begin{displaymath}
                     (Qx+C) \cdot (Qx+C) 
                      + \sigma x \cdot (\fourth \sigma x - C) 
                      \ < \ 1 \, 
                    \end{displaymath}
                    with $\sigma = \pm 4 \sqrt{ \, |K| \, }$${\bf ;}$  
              \smallskip 
              \item[$(-)_k$] 
                    or $h$ is the Klein model of sectional curvature $K$ on 
                    the unit ball in $\mathbb R^n$, 
                    and $W = Qx + C - ( x \cdot C ) x$ satisfies 
                    \begin{displaymath}
                     \frac{1}{\ 1 - (x \cdot x) \ } 
                     \{ (Qx+C) \cdot (Qx+C) - (x \cdot C)^2 \} 
                     \ < \ |K| \, .
                    \end{displaymath} 
             \end{itemize}
       \end{itemize}
       Furthermore, if $M$ is simply-connected and $h$ is complete, then 
       the said local isometry is in fact a global isometry. 
  \end{theorem}

  {\it Proof:}

  \begin{itemize} 
   \item By Proposition 1, every strongly convex Randers metric has the 
         representation, stipulated in $(1)$, in terms of the Zermelo 
         navigation variables $(h, W)$. 
   \smallskip 
   \item Theorem 3 tells us that $h$ must be a Riemannian space form. The 
         discussion after the statement of Theorem 3 reduces the landscape to 
         only four families, in keeping with (2). They are as follows.   
         \begin{itemize} 
          \item[$(+)$] 
               For $K>0$: $h$ must have sectional curvature $K$ and $W$ is 
               Killing. 
          \item[$(0)$] 
               For $K=0$: $h$ must be flat and $W$ is Killing. 
          \item[$(-)$] 
               For $K<0$: there are two scenarios,   
               \begin{itemize} 
                \item[$(-)_e$] 
                     either $h$ is flat, $\sigma = \pm 4 \sqrt{|K|}$, 
                     and $\mathcal L_W h = - \sigma \, h$ (in which case 
                     $W$ turns out to be $- \half \sigma$ times the radial 
                     vector $x = (x^i)$, plus an arbitrary Killing field);  
                \item[$(-)_k$] 
                     or $h$ has sectional curvature $K$ and $W$ is Killing. 
               \end{itemize} 
         \end{itemize} 
   \smallskip 
   \item Up to (Riemannian) isometry, there are only three standard 
         models for Riemannian metrics $h$ of constant sectional curvature 
         $K$.  They are: $\frac{1}{K}$ times the standard metric on the unit 
         $n$-sphere, Euclidean $\mathbb R^n$, and the Klein metric with 
         sectional curvature $K$ on the unit ball in $\mathbb R^n$. In view of 
         Lemma 2, when classifying $F$ up to Finslerian isometry, it suffices 
         to list the allowable vector fields $W$ with respect to each of the 
         three specific models.  (A more leisurely discussion of this point 
         is given at the beginning of \S 7.)  
         For the families $(+)$ and $(-)_k$, this has been done by 
         Proposition 6. Families $(0)$ and $(-)_e$ are handled by 
         Proposition 5, with $\sigma = 0$ and $\sigma = \pm 4 \sqrt{|K|}$, 
         respectively. 
   \smallskip 
   \item In each of the four families, the constraint that must be satisfied 
         by $Q$, $C$ and $x$ is equivalent to $|W| < 1$, which characterizes 
         the strong convexity of the Randers metric in question.  The table in 
         \S 5.4 shows that this constraint admits non-trivial solutions for 
         all four families.  In \S 7.2--7.4 we enumerate, with the help of 
         normal forms, all the $Q$, $C$ for which there exists an open domain 
         of $x$ on which $|W|<1$ holds. 
   \smallskip 
   \item Finally, if $M$ is simply-connected and $h$ is complete, Hopf's 
         classification theorem assures us that the Riemannian space form 
         $(M,h)$ must be globally isometric to one of the three standard 
         models. \qed 
  \end{itemize} 

                  %-----------------------------------------

 \subsection{Globally defined solutions on the standard $S^n$}  

  We see in the previous section that all strongly convex Randers metrics of 
  constant flag curvature $K > 0$ arise locally as solutions to Zermelo's 
  problem of navigation on the unit sphere $S^n$, under the influence of a 
  Killing field (an infinitesimal isometry) of $\frac{1}{K}$ times the 
  standard metric on $S^n$.  Let us show that each strongly convex solution on 
  any {\it closed} hemisphere has a unique smooth extension to a globally 
  defined strongly convex solution on $S^n$.  There is no restriction on the 
  dimension $n$. 
  
                  %.........................................

  \subsubsection{An extension} 

   Without loss of generality, let us assume that the hemisphere in question 
   is the {\it closed} eastern hemisphere.  Parametrize the eastern ($s = +1$) 
   and western ($s = -1$) open hemispheres, as submanifolds of the ambient 
   $\mathbb R^{n+1}$, by the maps     
   \begin{displaymath}
    x \mapsto \psi^{\pm}(x) := \frac{1}{\sqrt{1 + x \cdot x}} ( s , x ) \, , 
    \ \hbox{with} \ x \in \mathbb R^n \, .
   \end{displaymath} 
   Geometrically, the tangent space at the east pole (resp. west pole) is 
   identified with $\mathbb{R}^n$.  Each point $q$ on an open hemisphere lies 
   on a unique ray which emanates from the center of the sphere.  This ray 
   intersects the copy of $\mathbb R^n$ tangent to the pole, at a point $x$.  
   The above parametrization expresses $q$ in terms of $x$.  

   According to Theorem 7, on the {\it open} eastern hemisphere, the given 
   Randers metric has navigation data $(h,W)$, where $h$ is $\frac{1}{K}$ 
   times the standard Riemannian metric of $S^n$, and 
   $W(x) = Qx + C + (x \cdot C)x$.  We find that it is easier to visualize 
   $W(x)$ by considering its image under $\psi^+_*$.  Motivated by a 
   Lie-theoretic reason that will be pointed out in \S 7.2, we convert the 
   image point $p := \psi^+(x)$ into a position {\it row} vector $p^t$ of 
   $\mathbb R^{n+1}$.  A computation gives  
   \begin{displaymath}
    [\psi^+_*W(x)]^t 
    = 
    p^t \, \Omega \, , 
    \qquad \hbox{where} \qquad
    \Omega = 
    \left(
    \begin{array}{cc}
      0   &  C^t \\
     -C   &  -Q
    \end{array}
    \right) 
   \end{displaymath}
   is an $(n+1) \times (n+1)$ skew-symmetric {\it constant} matrix, and 
   ${}^t$ means transpose.  The continuity of $W$ on the closed hemisphere 
   implies that its value at any point $p$ on the equator is also the matrix 
   product $p^t \, \Omega$.  

   Extend $W$ to the open western hemisphere by insisting that 
   $[\psi^-_*W(x)]^t = [\psi^-(x)]^t \, \Omega$.  The result is 
   $W(x) = Qx + sC + (x \cdot sC)x$, with $s = -1$.  

   It is an artifact of local coordinates that $W$ is constructed from the 
   data $(Q,C)$ on the eastern hemisphere, but from $(Q,-C)$ on the western 
   hemisphere.  The actual Killing field on the embedded unit sphere in 
   $\mathbb R^{n+1}$ has the value $p^t \, \Omega$ at any point $p$, including 
   the equator.  Since the matrix $\Omega$ is constant, there is no question 
   that the constructed $W$ is globally defined and smooth. 

            %............................................
  
  \subsubsection{Uniqueness of the extension} 

   Let $W$ be any global extension of the given Killing field.  The 
   isometries of $(S^n,h)$ consist of rigid rotations, implemented by constant 
   $(n+1) \times (n+1)$ orthogonal matrices right multiplying the row vectors 
   of $\mathbb R^{n+1}$.  Since $W$ is an infinitesimal isometry, it is the 
   initial tangent to a curve of isometries.  Thus it also corresponds to a 
   {\it constant} matrix which right multiplies {\it all} row vectors.  For 
   points $p$ of the eastern hemisphere, we have determined the matrix in 
   question to be the above $\Omega$.  Constancy dictates that the same 
   $\Omega$ must be used for the western hemisphere as well.  This proves that 
   every global extension agrees with the one we presented.  In particular, 
   any global $W$ with data $(Q,C)$ on some hemisphere must have data $(Q,-C)$ 
   on the complement. 

            %............................................
   
  \subsubsection{Strong convexity} 

   The strong convexity criterion reads $|W| < 1$.  On the two {\it open} 
   hemispheres, Proposition 6 helps us deduce that 
   \begin{displaymath}
    |W(x)|^2 = 
    \frac{1}{K \{ \ 1 + (x \cdot x) \} \ }
    \{ (Qx+sC) \cdot (Qx+sC) + (x \cdot sC)^2 \} \, . 
   \end{displaymath}
   Using this formula, it is straight forward to check that 
   $|W(x)|^2 = (p^t \, \Omega) \cdot (p^t \, \Omega)$, where 
   $p = \psi^{\pm}(x)$.  Before the 
   extension, our Randers metric is strongly convex on the closed eastern 
   hemisphere.  In particular, $(p^t \, \Omega) \cdot (p^t \, \Omega) < 1$ for 
   all points $p$ of the open eastern hemisphere.  Replacing $p$ by $-p$ 
   generates all the points of the open western hemisphere, but does not alter 
   $(p^t \, \Omega) \cdot (p^t \, \Omega)$.  Therefore the extended metric is 
   also strongly convex on the open western hemisphere and hence on all of 
   $S^n$. 

            %...........................................

  \subsubsection{Discussion} 

   The examples of \S 3.1.1 and \S 3.1.2 determine globally defined Randers 
   metrics of constant positive flag curvature on $S^3$.  The first example 
   illustrates the necessity of assuming strong convexity on a closed 
   hemisphere.  Had we permitted $\tau = 1$, the norm of $W$ would have been 
   less than 1 on the open (eastern and western) hemispheres; but strong 
   convexity would fail at the points $(0,p_1,p_2,p_3)$ on the equator.

            %...........................................

 \subsection{Globally defined solutions on Euclidean $\mathbb{R}^n$}

  Because Euclidean $\mathbb{R}^n$ is covered by a single coordinate chart, 
  globality is relatively easy to address.  According to scenarios $(0)$ and 
  $(-)_e$ of Theorem 7, navigation on $\mathbb{R}^n$ under an infinitesimal 
  homothety $W$ produces a strongly convex Randers metric of constant flag 
  curvature $K \leqslant 0$ wherever $|W|<1$.  In particular, the Randers 
  metric is defined globally if and only if 
  \begin{displaymath}
    |W(x)|^2 
    = (Qx + C) \cdot (Qx + C) + \sigma x \cdot ( \fourth \sigma x - C )
    < 1 
    \quad 
    \hbox{for all} \ 
    x \in \mathbb{R}^n \, . 
  \end{displaymath}
  Here, $\sigma$ is zero if $K=0$, and has the values $\pm 4 \sqrt{|K|}$ if 
  $K<0$.  Since $|W(x)|^2$ is polynomial in $x$, the displayed criterion is 
  possible if and only if both $\sigma$ and $Q$ vanish, in which case $W = C$, 
  with $C \cdot C < 1$.  The resulting Randers metric is locally Minkowski.

  This conclusion is consistent with \S 3.2, where the only globally 
  defined example is that of \S 3.2.3.

            %...........................................

 \subsection{Globally defined solutions on the Klein model}  

  It remains to discuss global solutions to Zermelo's problem of navigation on 
  the Klein model with constant sectional curvature $K<0$, under the influence 
  of a Killing vector field $W$.  Theorem 7 says that the resulting Randers 
  metric has constant negative flag curvature $K$.  Strong convexity of the 
  Randers metric is equivalent to $|W|<1$.  In this subsection we will show 
  that requiring strong convexity on the entire open unit ball forces $W = 0$, 
  whence the negatively curved Randers metric is simply the Klein model 
  itself. 

  Suppose $|W|<1$ holds on the entire open unit ball.  It is implicit in 
  Proposition 6 that  
  \begin{displaymath}
   |W(x)|^2 = \frac{(Qx + C) \cdot (Qx + C) - (x \cdot C)^2}
                   {|K| \, ( 1 - x \cdot x )} \, .
  \end{displaymath}
  Note that $|K| ( 1 - x \cdot x ) > 0$ because $K$ is negative and $x$ is 
  confined to the unit ball.  Multiplying the inequality 
  $0 \leqslant |W|^2 < 1$ by this positive denominator yields 
  \begin{displaymath}
   0 
   \leqslant 
   (Qx + C) \cdot (Qx + C) - (x \cdot C)^2  <  |K| \, ( 1 - x \cdot x ) \, .
  \end{displaymath}
  Letting $x \cdot x \to 1$ leads to  
  $(Qx + C) \cdot (Qx + C) - ( x \cdot C )^2 = 0$ for all unit $x$.  In 
  particular, $(Qx + C) \cdot (Qx + C) = (-Qx + C) \cdot (-Qx + C)$, which is 
  equivalent to $Qx \cdot C = 0$.  The equality above then simplifies to  
  $Qx \cdot Qx + C \cdot C - ( x \cdot C )^2 = 0$, again for all unit $x$.  

  Since we are in dimension at least two, there exists a unit $x_0$ such that 
  $x_0 \cdot C = 0$.  The ensuing equation 
  $Qx_0 \cdot Qx_0 + C \cdot C = 0$ tells us that $C$ must have been zero to 
  begin with.  This reduces our original equality to $Qx = 0$ for all unit 
  $x$, implying that $Q = 0$.  Thus $W$ is identically zero, and our assertion 
  follows. 
  
%------------------------------------------------------------------------------

\section{The moduli space} 

           %-----------------------------------------------

 \subsection{A strategy} 

  Theorem 3 of Section 4.4 characterizes the navigation data $(h,W)$ of 
  strongly convex Randers metrics with constant flag curvature $K$.  Namely, 
  $h$ must be a Riemannian metric with constant sectional curvature 
  $K + \onesixteenth \sigma^2$, and $W$ must be an infinitesimal homothety 
  of $h$.  Also, we observed that $\sigma$ can be nonzero only when $h$ is 
  flat. 

  Consider any Randers metric $(M,F)$ of constant flag curvature $K$, with 
  navigation data $(h,W)$.  There exists a local isometry $\varphi$ between 
  $(M,h)$ and one of the three standard models: 
  \begin{itemize}
   \item the sphere $(S^n,h_+)$ of constant curvature $K$ when $K>0$;
   \item Euclidean space $(\mathbb{R}^n, h_0)$ when $K = 0$, 
         or when $K<0$ and $\sigma = \pm 4 \sqrt{|K|}$; 
   \item the Klein model $(\mathbb{B}^n,h_-)$ of constant curvature $K$ 
         when $K<0$ and $\sigma = 0$.
  \end{itemize}
  By using this $\varphi$ to transform the navigation data $(h,W)$ if 
  necessary, we may assume without loss of generality that $h$ is already one 
  of the standard models.  For each such $h$, Theorem 7 of Section 6.1 lists 
  its infinitesimal homotheties $W$.  

  That list contains a good amount of redundancy because it includes 
  Randers metrics that are locally isometric.  The redundancy comes from the 
  symmetry/isometry group $G$ of $h$, consisting of diffeomorphisms $\phi$ 
  that leave $h$ invariant.  Since $\phi^* h = h$, the action of the Lie group 
  $G$ on the navigation data is $(h,W) \mapsto (h,\phi_*W)$.  According to 
  Lemma 2 of Section 2.4, sets of navigation data which lie on the same 
  $G$-orbit correspond to locally isometric Randers metrics.  The redundancy 
  we described can therefore be eliminated by collapsing each $G$-orbit to a 
  point.  These ``points'' constitute the elements of our moduli space 
  $\mathcal M_K$ of Randers metrics with constant flag curvature $K$.  It is 
  the goal of \S 7 to parametrize $\mathcal M_K$ and thereby count 
  its dimension. 

  To this end, we begin with a standard model $h$ ($=$ $h_+$, or $h_0$, or 
  $h_-$) of a Riemannian space form.  Identify the isometry group $G$ of $h$ 
  with a matrix subgroup of $GL_{n+1}\mathbb{R}$.  The infinitesimal 
  homotheties $W$ of $h$ comprise a representation of some matrix Lie 
  subalgebra $\mathfrak{h}$ of $\mathfrak{gl}_{n+1}\mathbb{R}$.  The 
  push-forward action $W \mapsto \phi_* W := \phi_* \circ W \circ \phi^{-1}$ 
  then corresponds to the ``adjoint action"   
  \begin{displaymath}
   \Omega \ \mapsto \ \hbox{Ad}_{g} \Omega := g \, \Omega \, g^{-1} 
  \end{displaymath} 
  of $G$ on $\mathfrak h$.  Here: 
  {\bf (1)} 
      $g \in GL_{n+1}\mathbb{R}$ is the matrix which corresponds to the 
      isometry map $\phi$, and $\Omega \in \mathfrak{h}$ is the matrix analog 
      of the infinitesimal homothety $W$ (which is a vector field).  
  {\bf (2)} 
      $Ad$ is well defined because the equation $\mathcal L_W h = - \sigma h$, 
      being tensorial, becomes $\mathcal L_{\phi_*W} h = - \sigma h$.  Thus, 
      $\phi_*W$ is an infinitesimal homothety of $h$ whenever $W$ is, and the 
      value of $\sigma$ is invariant under isometries. 
  {\bf (3)} 
      According to Theorem 3, when $h$ is not flat, its infinitesimal 
      homotheties are simply its Killing vector fields.  In that case, 
      $\mathfrak{h}$ equals the Lie algebra $\mathfrak g$ of $G$, and 
      $Ad$ is the standard adjoint action of a Lie group on its Lie algebra. 
 
  The adjoint action $Ad$ described above partitions $\mathfrak h$ into 
  orbits.  These orbits correspond to distinct local isometry classes of 
  Randers metrics with constant flag curvature $K$.  For each orbit, 
  matrix theory singles out a privileged representative $\tilde \Omega$, to 
  be referred to as a normal form.  These normal forms provide a concrete 
  parametrization of the points in the moduli space $\mathcal M_K$, and the 
  number of parameters constitutes its dimension.  The linear algebra behind 
  the construction of $\mathcal M_K$ depends on the sign of $K$.  Here is an 
  overview.  
  \begin{itemize} 
   \item For $K > 0$, $h = h_+$ is $\frac{1}{K}$ times the standard metric on 
         the unit $n$-sphere.  The orbits are those which result from the 
         adjoint action of the orthogonal group $O(n+1)$ on its Lie algebra 
         $\mathfrak o(n+1)$. 
   \item For $K = 0$, we have $h = h_0$, the standard flat metric on 
         $\mathbb R^n$.  The orbits come from the adjoint action of the 
         Euclidean group $E(n)$ on its Lie algebra $\mathfrak e(n)$.  Here, 
         $E(n)$ is the semi-direct product of $O(n)$ with the additive group 
         $\mathbb R^n$ of translations. 
   \item For $K < 0$, the orbits consist of two camps.  (1) $h = h_-$ is the 
         Klein model, and the $Ad$ orbits arise from a subgroup of the Lorentz 
         group $O(1,n)$, acting on the Lie algebra $\mathfrak o(1,n)$.  
         (2) $h = h_0$ is the flat Euclidean metric, and the $Ad$ orbits are 
         those of $E(n)$ acting on a matrix description of the 
         infinitesimal homotheties, with $\sigma = \pm 4 \sqrt{|K|}$.   
  \end{itemize}
  The Lie theory necessary for determining the normal form $\tilde \Omega$ is 
  relegated to the Appendix (\S 10).  The material there will be called upon 
  frequently in the following three subsections as we determine the local 
  isometry classes of strongly convex Randers metrics of constant flag 
  curvature.

                %-----------------------------------------

 \subsection{The $n$-sphere}

  The isometry group $G$ of $(S^n,h_+)$ is $O(n+1)$, whose elements are 
  orthogonal matrices which implement rigid rotations by right multiplying the 
  row vectors of $\mathbb R^{n+1}$.  As explained in \S 6.2, each Killing 
  vector field $W$ of $(S^n,h_+)$ also corresponds to a constant matrix which 
  right multiplies those row vectors, and we have identified that 
  skew-symmetric $(n+1) \times (n+1)$ matrix to be 
  \begin{displaymath}
   \Omega := 
   \left(
   \begin{array}{cc}
    0  & C^t \\
    -C & -Q 
   \end{array}
   \right) \, ,   
  \end{displaymath}
  an element of the Lie algebra $\mathfrak o(n+1)$.  This correspondence 
  between the Killing fields of $(S^n,h_+)$ and $\mathfrak o(n+1)$ is a Lie 
  algebra isomorphism.  (Incidentally, if we had let the group $O(n+1)$ act on 
  column vectors instead, then the matrix $- \Omega$ would correspond to $W$,  
  while the {\it negative} of the commutator $[-\Omega_1, -\Omega_2]$ would 
  represent the Lie bracket $[W_1,W_2]$, rendering the correspondence a Lie 
  algebra anti-isomorphism.)
  
  Applying \S 10.2 (with $\ell := n+1$) to $\Omega$, we see that there exists 
  a $g \in O(n+1)$ so that $\tilde \Omega = g \Omega g^{-1}$ is in normal 
  form.  Explicitly: 
  \begin{center} 
   \begin{tabular}{lll} 
    when $n$ is even, 
    & 
    $\tilde \Omega = a_1 J \oplus \cdots \oplus a_m J \oplus 0$
    & 
    with $m = n / 2$; \\ 
    when $n$ is odd, 
    &
    $\tilde \Omega = a_1 J \oplus \cdots \oplus a_m J$
    &
    with $m = (n+1)/2$.  
   \end{tabular}  
  \end{center}  
  Here, $a_1 \geqslant a_2 \geqslant \cdots \geqslant a_m \geqslant 0$ and 
  \begin{displaymath}
   J = \left(
   \begin{array}{rr}
   0  & 1 \\
   -1 & 0 \\ 
   \end{array}
   \right) \, .
  \end{displaymath}    

  The matrix $\tilde \Omega$ represents the Killing field 
  $\tilde W = \phi_* W$, where $\phi$ is the map which corresponds to the 
  orthogonal matrix $g$ (\S7.1).  According to Theorem 7, 
  $\tilde W$ has the form $\tilde Q x + \tilde C + (x \cdot \tilde C) x$ with 
  respect to the projective coordinates $x$ which parametrize the eastern 
  hemisphere.  Comparing the matrix analog 
  \begin{displaymath}
   \left(
    \begin{array}{cc}
     0          &  \tilde C^t \\
     -\tilde C  &  -\tilde Q
    \end{array}
   \right) 
  \end{displaymath}
  of $\tilde W$ with $\tilde \Omega$, we conclude that 
  $\tilde C^t  = ( a_1 , 0, \dots, 0)$ and  
  $-\tilde Q = 0 \oplus a_2 J \oplus \cdots \oplus a_m J \oplus 0$ when 
  $n$ is even, $-\tilde Q = 0 \oplus a_2 J \oplus \cdots \oplus a_m J$  
  when $n$ is odd.  

  The Randers metric which solves Zermelo's problem of navigation on 
  $(S^n,h_+)$ under the influence of $\tilde W$ must satisfy the strong 
  convexity criterion $|\tilde W|<1$.  In terms of the data 
  $(\tilde Q, \tilde C)$ for $\tilde W$, inequality (2,+) of Theorem 7 
  expresses this criterion as:   
  \begin{center} 
   \begin{tabular}{ll}
    \quad 
    $a_1^2 ( 1 + x_1^2 ) + a_2^2 ( x_2^2 + x_3^2 ) + \cdots + 
     a_m^2 ( x_{n-2}^2 + x_{n-1}^2 ) < K ( 1 + x \cdot x ),$
    &
    $n$ even; \\
    \quad  
    $a_1^2 ( 1 + x_1^2 ) + a_2^2 ( x_2^2 + x_3^2 ) + \cdots + 
     a_m^2 ( x_{n-1}^2 + x_n^2 ) \quad < K ( 1 + x \cdot x ),$     
    &
    $n$ odd.
   \end{tabular} 
  \end{center}  
  We wish to demarcate those $a_i$ that satisfy the above inequalities on an 
  open set.

           %----------------------------------------

   \subsubsection{Locally defined metrics when $n$ is even}

    Consider the point $x = (0, \dots, 0, x_n)$.  Here the condition 
    $|\tilde W(x)|<1$ simplifies to $a_1^2 < K (1 + x_n^2)$, which can be 
    made to hold for arbitrary but fixed $a_1$ by choosing $|x_n|$ large 
    enough.  Once we have $|\tilde W (x)| < 1$, the continuity of $\tilde W$ 
    effects $|\tilde W|<1$ on a neighborhood about this $x$.  Thus, for even 
    $n$, the moduli space is parametrized by 
    \begin{displaymath}
     a_1 \geqslant \ldots \geqslant a_m \geqslant 0 
    \end{displaymath}
    without any upper bound on $a_1$, and hence none on the $a_i$.     
 
     %------------------------------------------------------
    
   \subsubsection{Locally defined metrics when $n$ is odd}

    Suppose $|W|<1$ holds at some point $x$.  Then 
    $0 \leqslant a_m \leqslant a_i$ implies that 
    \begin{center}
     \begin{tabular}{l}
      $a_m^2 ( 1 + x \cdot x ) 
       \leqslant  
       a_1^2 (1 + x_1^2) + a_2^2 (x_2^2 + x_3^2) + \cdots + 
       a_m^2 (x_{n-1}^2 + x_n^2)                                  
       <  
       K (1 + x \cdot x)$.
     \end{tabular} 
    \end{center}
    In particular, we obtain the necessary condition $a_m < \sqrt{K}$.  
    Conversely, given $a_m < \sqrt{K}$, let us consider a point $x$ of the 
    form $(0, \dots, 0, x_n)$.  At this $x$, the desired condition 
    $|\tilde W(x)|<1$ simplifies and can be rearranged to read 
    $a_1^2 < K + (K - a_m^2) x_n^2$.  Since $a_m^2 < K$, the inequality can be 
    made to hold by choosing $|x_n|$ large enough.  Continuity then extends 
    $|\tilde W|<1$ from this $x$ to a neighborhood containing it.  Therefore 
    the isometry classes of locally defined Randers metrics on the odd 
    dimensional spheres are parametrized by 
    \begin{displaymath}
     a_1 \geqslant \ldots \geqslant a_m \geqslant 0 \, , 
     \quad \textrm{with} \quad a_m < \sqrt{K} \, .
    \end{displaymath}

     %------------------------------------------------------

   \subsubsection{Globally defined metrics}

    Here, the criterion $|\tilde W(x)| < 1$ must hold on the entire sphere.  
    In particular, it must hold for all $x \in \mathbb R^n$ parametrizing 
    the open eastern hemisphere.  Setting $x=0$ in the inequalities 
    immediately before \S 7.2.1 gives $a_1 < \sqrt{K}$.  Conversely, if 
    $a_1 < \sqrt{K}$, then those inequalities are satisfied for all $x$ 
    because $a_1 \geqslant a_i \geqslant 0$.  
    Hence the constraint $a_1 < \sqrt{K}$ is both necessary and sufficient 
    for strong convexity on the open eastern hemisphere.  By virtue of  
    \S 6.2.3, the same bound on $a_1$ effects $|\tilde W|<1$ on the 
    open western hemisphere.  Thus strong convexity holds on the open 
    hemispheres if and only if the condition $a_1 < \sqrt{K}$ is met.   

    It turns out that $a_1 < \sqrt{K}$ ensures strong convexity on the equator 
    as well.  To see this, let $u$ be any unit vector in the copy of 
    $\mathbb R^n$ tangent to the poles.  Our parametrization (see \S 6.2.1) of 
    the open hemispheres says that $\lim_{t \to \infty} tu$ corresponds 
    asymptotically to some point $p$ on the equator.  In fact,  
    $p 
     = \lim_{t \to \infty} (1 + tu \cdot tu)^{-1/2} (s,tu)  
     = (0,u)$.  
    Calculating with the norm $|y|^2 := h(y,y)$ given in part (a) of 
    Proposition 6, we find that  
    \begin{displaymath} 
     |\tilde W(p)| 
     = \lim_{t \to \infty} |\tilde W(tu)| 
     = \frac{1}{\sqrt{K}} \sqrt{ ( \, u \cdot s \tilde C \, )^2 + 
                                 | \, \tilde Qu \, |^2 }\, , 
    \end{displaymath} 
    which is independent of $s = \pm 1$.  A direct computation, 
    using the fact that $a_1$ dominates all other $a_i$, and 
    $u \cdot u = 1$, yields $|\tilde W(p)|^2 \leqslant (a_1)^2 / K < 1$.   

    Thus the moduli space for the isometry classes of globally defined 
    constant flag curvature $K > 0$ Randers metrics on $S^n$ is given by the 
    polytope
    \begin{displaymath}
     \sqrt{K} > a_1 \geqslant \cdots \geqslant a_m \geqslant 0 \, .
    \end{displaymath}

   \subsubsection{Global versus local}

    For the locally defined metrics, the upper bound $a_1 < \sqrt{K}$ 
    is not necessary because the strong convexity criterion $|\tilde W| < 1$ 
    only has to hold on some open subset of $S^n$.  However, when $n$ is odd, 
    all local solutions have to satisfy $a_m < \sqrt{K}$. 

    The metric of \S 3.1.1 illustrates these nuances well.  The table in 
    \S 5.4 tells us that $C^t = (0,0,0)$ and $Q = \tau J \oplus 0$.  Using 
    the data $(Q,C)$, construct $\Omega$ as in \S 7.2.  Almost by inspection, 
    the normal form is $\tilde \Omega = \tau J \oplus 0 J$, thus $a_1 = \tau$ 
    and $a_m \equiv a_2 = 0$.  Since $K$ here is $1$, the theory assures us 
    that a locally defined strongly convex solution exists for any $\tau$, 
    while strongly convex global solutions are characterized by $\tau < 1$. 

    Indeed, \S 6.2.3 tells us that $\tilde W(p) = p^t \, \tilde \Omega$, and 
    $|\tilde W(p)|^2 = (p^t \, \tilde \Omega) \cdot (p^t \, \tilde \Omega) 
    = \tau^2 (p_0^2+p_1^2)$, where $p^t = (p_0,p_1,p_2,p_3)$ gives the  
    coordinates of an {\it arbitrary} point on the embedded $S^3$ in 
    $\mathbb R^4$.  So $|\tilde W|<1$ globally, as long as $\tau < 1$.  On 
    the other hand, if $\tau \geqslant 1$, then $|\tilde W(p)|<1$ holds only 
    at those points $p$ on $S^3$ where $p_0^2 + p_1^2 < 1 / \tau^2$.

                   %..................................

   \subsubsection{The moduli space for K positive}

    \begin{proposition}
     The local isometry moduli space of $n$-dimensional strongly convex 
     Randers spaces of constant flag curvature $K > 0$ is parametrized by 
     $a = (a_1 , \ldots , a_m ) \in \mathbb{R}^m$ as follows.
     \begin{itemize}
      \item[$\circ$] When $n$ is even, $m=n/2$ and the parameter space is 
                     given by 
                     \begin{displaymath}
                      a_1 \geqslant \cdots \geqslant a_m \geqslant 0 \, .
                     \end{displaymath}
      \item[$\circ$] When $n$ is odd, $m=(n+1)/2$ and the parameter space 
                     is given by 
                     \begin{displaymath}
                      a_1 \geqslant \cdots \geqslant a_m \geqslant 0 \ , 
                      \quad \textrm{with} \ \sqrt{K} > a_m \, .
                     \end{displaymath}
      \item[$\circ$] The globally defined metrics on $S^n$ are parametrized 
                     by the polytope 
                     \begin{displaymath}
                      \sqrt{K} > a_1 \geqslant \cdots \geqslant a_m 
                                                      \geqslant 0 \, .
                     \end{displaymath}
     \end{itemize}
    \end{proposition}

                  %--------------------------------------

 \subsection{Euclidean space}

  The isometry group of $( \mathbb{R}^n , h_0 )$ consists of rotations, 
  reflections, and translations; it is the Euclidean group $E(n)$.  Though the 
  action of $E(n)$ on $\mathbb R^n$ is affine, it can be implemented by matrix 
  multiplication.  To this end, we first represent elements $\phi$ of $E(n)$ 
  by matrices $g \in GL_n\mathbb{R}$ of the form
  \begin{displaymath} 
   g = 
   \left( 
    \begin{array}{cc}
     A & 0 \\
     b & 1
    \end{array}
   \right) \, , 
   \quad \hbox{where} \quad
   A \in O(n) \hbox{ and } b \in \mathbb{R}^n \, .
  \end{displaymath} 
  Next we embedd Euclidean $n$-space into $\mathbb R^{n+1}$ by assigning to 
  each point $x$ the column position vector $\psi(x) = (x,1) =: p$.  The 
  matrix action we have in mind is then 
  \begin{displaymath}
    p^t \mapsto p^t g = (x^t A + b,1) \, .
  \end{displaymath}
  Here, $p^t$ and the output $p^t g$ are both row vectors. 

  The image of an infinitesimal homothety $W = - \half \sigma x + Qx + C$ 
  under the described representation is $[\psi_* W(x)]^t = p^t \, \Omega$, 
  where
  \begin{displaymath}  
   \Omega := \left( 
             \begin{array}{cc}
              -\half \sigma I_n - Q & 0 \\
                          C^t       & 0
             \end{array}
             \right) \, .
  \end{displaymath} 
  Such matrices, with $\sigma \in \mathbb{R}$, $C \in \mathbb{R}^n$ and 
  $Q \in \mathfrak o(n)$, form a Lie subalgebra $\mathfrak h$ of 
  $\mathfrak{gl}_{n+1}$.  The correspondence between the infinitesimal 
  homotheties $W$ of $( \mathbb{R}^n , h_0 )$ and the subalgebra 
  $\mathfrak h$ is a Lie algebra isomorphism.  When $\sigma=0$, $\mathfrak h$ 
  is the Lie algebra $\mathfrak e(n)$ of $E(n)$. 

  The vector field 
  $\tilde W = - \half \tilde \sigma x + \tilde Q x + \tilde C$ is the push 
  forward of $W$ under an isometry $\phi \in G$ if and only if its matrix 
  representative $\tilde \Omega$ is given by $g \Omega g^{-1}$.  Since 
  \begin{displaymath} 
   g^{-1} = 
   \left(
    \begin{array}{cc}
     A^t    & 0 \\
     -bA^t  & 1
    \end{array}
   \right) \, ,
  \end{displaymath} 
  we have 
  \begin{displaymath} 
   \left(
    \begin{array}{cc}
     - \half \tilde \sigma I_n - \tilde Q & 0 \\
             \tilde C^t                   & 0
    \end{array}
   \right)  
   =
   \tilde \Omega 
   =    
   g \Omega g^{-1} 
   = 
   \left(
    \begin{array}{cc}
     - \half \sigma I_n - A Q A^t & 0 \\
       \left[ A W(b) \right]^t    & 0 
    \end{array}
   \right) \, ,
  \end{displaymath} 
  where $W(b) = - \half \sigma b + Qb + C$.  Thus $\tilde \sigma = \sigma$, 
  $\tilde Q = A Q A^t$, and $\tilde C = A W(b)$; in particular, the value of 
  $\sigma$ remains unchanged under any isometry, a general fact we pointed out 
  in \S 7.1.  Our objective is to find $A$ and $b$, equivalently $g \in E(n)$, 
  so that $\tilde \Omega$ takes on a simplest form.    

               %.......................................

 \subsubsection{The case of $\sigma = 0$ and the moduli space for $K=0$}  

  The Randers metrics of constant flag curvature zero arise as perturbation of 
  the Euclidean metric under an infinitesimal isometry.  This corresponds to 
  the $\sigma=0$ case in the above discussion.  

  For ease of exposition, let us abbreviate group elements $g \in E(n)$ as 
  $\{ A, b \}$ and Lie algebra elements $\Omega \in \mathfrak e(n)$ as 
  $[ -Q, C^t ]$.  
  \begin{itemize} 
   \item[(1)] By \S 10.2, we can find an $R \in O(n)$ which puts $-Q$ into 
              the normal form $- \tilde Q = 
              \rho_1 J \oplus \cdots \oplus \rho_h J \oplus 0_{n-2h}$, with 
              $\rho_1 \geqslant \cdots \geqslant \rho_h > 0$.  Thus  
              $g_1 := \{ R, 0 \}$ conjugates $\Omega$ into 
              $\tilde \Omega_1 := [ -\tilde Q, (RC)^t ]$. 
   \smallskip 
   \item[(2)] Choose $r \in O(n-2h)$ to transform the last $n-2h$ components 
              of $RC$ into $(0, \dots, 0, \xi \geqslant 0)$, without affecting 
              its first $2h$ components $D := (D_1, \dots, D_h)$, listed 
              pairwise for convenience as $D_i = [C_{2i-1},C_{2i}]$.  The 
              corresponding group element $g_2 := \{ I_{2h} \oplus r, 0 \}$ 
              conjugates $\tilde \Omega_1$ into 
              $\tilde \Omega_2 := [-\tilde Q, (D, 0, \dots, 0, \xi)^t]$. 
   \smallskip 
   \item[(3)] Set $b := ( \frac{-JD_1}{\rho_1}, \dots, \frac{-JD_h}{\rho_h}, 
                          0, \dots, 0)$ and note that $-\tilde Q b = 
              (D, 0, \dots, 0)$.  Then $g_3 := \{ I_n, b \}$ conjugates 
              $\tilde \Omega_2$ into $\tilde \Omega_3 := 
              [-\tilde Q, (0, \dots, 0, \xi)^t]$.  
  \end{itemize} 
  In short, using $g := g_3 g_2 g_1 \in E(n)$, we get
  \begin{displaymath} 
   \tilde \Omega 
   := g \Omega g^{-1} 
    = 
    \left(
     \begin{array}{ccc}
      \rho_1 J \oplus \dots \oplus \rho_h J  &  0                & 0   \\
      0                                      &  0_{n-2h}         & 0   \\ 
      0                                      &  0, \dots, 0, \xi & 0 
     \end{array}
    \right) \, . 
  \end{displaymath} 
  A moment's thought tells us that $\xi = 0$ whenever 
  $C \in \textrm{Range} \, Q$, and $\xi > 0$ otherwise. 
  The strong convexity condition $|\tilde W|<1$ restricts our domain to 
  those $x$ which satisfy  
  \begin{displaymath}
   |\tilde W(x)|^2 
   =  ( \tilde Q x + \tilde C ) \cdot ( \tilde Q x + \tilde C ) 
   = \xi^2 + \sum_{i=1}^h \rho_i^2 (x_{2i-1}^2 + x_{2i}^2) < 1 \, .       
  \end{displaymath}
  In particular, we must have $\xi < 1$.  As long as this condition is met, 
  the inequality above will always be satisfied on some neighborhood of 
  the origin in $\mathbb R^n$.

  Suppressing the rank of $Q$ by augmenting the parameters, followed by some 
  appropriate relabeling, simplifies the normal form $\tilde \Omega$ to  
  \begin{displaymath} 
    \left(
     \begin{array}{cc}
      a_1 J \oplus \dots \oplus a_m J  &  0  \\
      0, \dots \dots, 0, a_0           &  0 
     \end{array}
    \right)_{\textrm{for even} \, n} \, , 
   \quad \quad 
    \left(
     \begin{array}{ccc}
      a_2 J \oplus \dots \oplus a_m J  &  0        & 0  \\
      0                                &  a_1      & 0 
     \end{array}
    \right)_{\textrm{for odd} \, n} \, . 
  \end{displaymath} 
  Here, {\it a priori} we have 
  \begin{center} 
   \begin{tabular}{llll} 
    $1 > a_0 \geqslant 0$, 
    &
    $a_1 \geqslant \cdots \geqslant a_m \geqslant 0$, 
    & 
    and $m = n / 2$ 
    &
    for even $n$;        \\ 
    $1 > a_1 \geqslant 0$, 
    &  
    $a_2 \geqslant \cdots \geqslant a_m \geqslant 0$, 
    &  
    and $m = (n+1) / 2$ 
    &    
    for odd $n$. 
   \end{tabular} 
  \end{center}  
  However: 
  \begin{itemize} 
   \item    When $n$ is even, $a_0$ and $a_m$ cannot both be nonzero for 
            any fixed $\Omega$.  Indeed, if $a_0 > 0$, then $C$ is not in  
            $\textrm{Range} \, Q$ and we must at least have $a_m = 0$.  
            On the other hand, if $a_m \not= 0$, then $Q$ is surjective and 
            thus $a_0$ must vanish.    
   \item    When $n$ is odd, the displayed normal form precludes any sort of 
            rigid coupling between $a_1$ and $a_m$. 
  \end{itemize} 
  For the even $n$ case, whenever $a_0 > 0$ (so that $a_m = 0$), let us agree 
  to relabel the remaining parameters $a_0, a_1, \dots, a_{m-1}$ as 
  $a_1, a_2, \dots, a_m$.  

  \begin{proposition}
   The local isometry moduli space of $n$-dimensional strongly convex Randers 
   spaces of constant flag curvature $K=0$ is parametrized by 
   $a = (a_1 , \ldots , a_m ) \in \mathbb{R}^m$ as follows.
   \begin{itemize}
    \item[$\circ$] 
          When $n$ is even, $m=n/2$ and the parameter space is 
          the disjoint union of 
          \begin{displaymath}
           a_1 \geqslant \cdots \geqslant a_m \geqslant 0          
           \qquad \hbox{and} \qquad 
           1 > a_1 > 0 \ , \ a_2 \geqslant \cdots a_m \geqslant 0 \, .
          \end{displaymath} 
    \item[$\circ$]  
          When $n$ is odd, $m=(n+1)/2$ and the parameter space is given by 
          \begin{displaymath}
           1 > a_1 \geqslant 0 \ , \ 
           a_2 \geqslant \cdots \geqslant a_m \geqslant 0 \, .
          \end{displaymath}
    \item[$\circ$] The globally defined metrics are parametrized by 
          \begin{displaymath}
           1 > a_1 \geqslant 0 \ , \ a_2 = \cdots = a_m = 0 \, .
          \end{displaymath}
   \end{itemize}
  \end{proposition} 
    
           %......................................

 \subsubsection{When $\sigma$ is nonzero} 

  Refer to the general discussion at the beginning of \S 7.3, and 
  the abbreviation introduced in \S 7.3.1.  We see that conjugating 
  $\Omega = [-\half \sigma I_n - Q, C^t]$ by any $g := \{ A, b \} \in E(n)$ 
  converts it to $[-\half \sigma I_n - A Q A^t, (A W(b))^t]$.  Select 
  $A \in O(n)$ to cast $-Q$ into the following normal form: 
  \begin{center} 
   \begin{tabular}{lll}  
    when $n$ is even, 
    &
    $-\tilde Q = - A Q A^t = a_1 J \oplus \cdots \oplus a_m J$, 
    & 
    with $m = n / 2$;      \\ 
    when $n$ is odd, 
    &
    $-\tilde Q = - A Q A^t = a_1 J \oplus \cdots \oplus a_m J \oplus 0$,
    &
    with $m = (n-1) / 2$. 
   \end{tabular} 
  \end{center} 
  Here, $a_1 \geqslant \cdots \geqslant a_m \geqslant 0$.  Note that   
  $W(b) = (Q - \half \sigma I_n)b + C$.  The linear operator 
  $Q - \half \sigma I_n$ is invertible because the spectrum of $Q$ is pure 
  imaginary (\S 10.2) whereas $\sigma$ is real and nonzero.  Therefore we may 
  select $b$ so that $W(b) = 0$.  With this choice of $A$ and $b$, 
  $g := \{ A, b \}$ conjugates $\Omega$ into the normal form
  \begin{displaymath} 
   \tilde \Omega 
   =    
   g \Omega g^{-1} 
   = 
   \left(
    \begin{array}{cc}
     - \half \sigma I_n - \tilde Q & 0 \\
     0                             & 0 
    \end{array}
   \right) \, . 
  \end{displaymath} 

  The corresponding infinitesimal homothety has $\tilde C = 0$ and its formula 
  is $\tilde W(x) = - \half \sigma x + \tilde Q x$.  Navigating on Euclidean 
  $\mathbb R^n$ subject to the wind $\tilde W$ generates a Randers metric of 
  negative flag curvature $K = -\onesixteenth \sigma^2$.  This metric is 
  strongly convex wherever 
  \begin{displaymath}
   |\tilde W(x)|^2 
   = 
   \tilde Q x \cdot \tilde Q x + \fourth \sigma^2 x \cdot x 
   = 
   \fourth \sigma^2 x \cdot x \ + \ 
   \sum_{i=1}^m a_i^2 (x_{2i-1}^2 + x_{2i}^2) \ < \ 1 \, .
  \end{displaymath}
  \begin{itemize} 
   \item[$\circ$] 
         For any choice of $a_i$ and $\sigma \not= 0$, this condition will be 
         satisfied on some neighborhood of the origin in $\mathbb{R}^n$.  
   \item[$\circ$] 
         The left-hand side is a nonzero polynomial in $x$.  Therefore strong 
         convexity will never hold globally on $\mathbb{R}^n$. 
   \item[$\circ$] 
         The space of local isometry equivalence classes is parametrized by 
         \begin{displaymath} 
          a_1 \geqslant \cdots \geqslant a_m \geqslant 0 \, ,
         \end{displaymath}
         with $m = n/2$ when $n$ is even, and $m = (n-1)/2$ when $n$ is odd.
  \end{itemize} 
  In order to complete our parametrization of the local isometry classes of 
  constant negative flag curvature Randers spaces, it remains to consider 
  perturbations of the Klein model.

                   %-----------------------------------

 \subsection{Hyperbolic space}

  In analogy with the spherical (\S 6.2, \S 7.2) and Euclidean (\S 7.3) cases, 
  we embed the Klein model of hyperbolic geometry into an ambient 
  $(n+1)$ dimensional space.  To that end, consider $\mathbb{R}^{n+1}$ 
  equipped with the scalar product $\langle v , w \rangle := v^t E w$, 
  where $E = -1 \oplus I_n$.  The isometry group of this space is the 
  Lorentz group $O(1,n)$.

  For $K<0$, define the subspace 
  $H_K := \{ x \in \mathbb{R}^{n+1} \ | \ 
             \langle x , x \rangle = \frac{1}{K} \}$.  
  We make three observations \cite{ON83}:  
  \begin{itemize}
   \item[$\circ$]  
         $H_K$ consists of two components, each diffeomorphic 
         to $\mathbb{R}^n$.
   \item[$\circ$]  
         $\langle \, , \, \rangle$ restricts to a Riemannian metric of 
         constant sectional curvature $K$ on $H_K$.
   \item[$\circ$]   
         $O(1,n)$ preserves $H_K$.
  \end{itemize}
  Let $h_K$ denote that component which passes through 
  $(1/\sqrt{|K|},0,\ldots,0)$.  Then $h_K$ is a complete, simply connected 
  model of hyperbolic space.  The isometry group $G$ of $h_K$ consists of 
  those matrices $g \in O(1,n)$ such that $g(h_K) = h_K$.  This identifies $G$ 
  as the orthochronous subgroup $O_+(1,n)$.  Its Lie algebra is 
  $\mathfrak{o}(1,n)$.

  Let us determine the relationship between Killing vector fields on the Klein 
  model and the Lie algebra $\mathfrak{o}(1,n)$.  Introduce the diffeomorphism
  \begin{displaymath}
   \psi (x) = \frac{1}{|K| \sqrt{ 1 - x \cdot x }} (1,x)
  \end{displaymath}  
  which maps the unit ball in $\mathbb{R}^n$ onto $h_K$.  The map $\psi$ is an 
  isometry between the Klein model and $h_K$.  Let $p := \psi(x)$ abbreviate 
  the position column vector of the image point.  Then  
  Killing vector fields $W(x) = Qx + C - ( x \cdot C ) x$ of the Klein model 
  are associated with elements 
  \begin{displaymath}
   \Omega := \left(
              \begin{array}{cc}
               0  &  C^t \\
               C  &  -Q  
              \end{array}
             \right)  
   \quad \in \ \mathfrak{o}(1,n)
  \end{displaymath}
  via $[\psi_* W(x)]^t = p^t \, \Omega$.  This correspondence is a Lie 
  algebra isomorphism.

  In \S 10.3 of the Appendix, we show that there exists a $g \in O_+(1,n)$ so 
  that $\tilde \Omega = g \Omega g^{-1}$ assumes one of three possible block 
  diagonal forms, as follows. 
  \begin{itemize}
   \item $i\Omega$ has a timelike eigenvector:
    \begin{center} 
    \begin{tabular}{lll} 
     when $n$ is even, 
     & 
     $\tilde \Omega = 0 \oplus a_1 J \oplus \cdots \oplus a_m J$, 
     & 
     with $m = n / 2$; \\ 
     when $n$ is odd, 
     &
     $\tilde \Omega = 0 \oplus a_1 J \oplus \cdots \oplus a_m J \oplus 0$, 
     &
     with $m = (n-1)/2$.  
    \end{tabular}  
   \end{center} 
   Here, $a_1 \geqslant a_2 \geqslant \cdots \geqslant a_m \geqslant 0$.
   See \S 3.3.1 for an example of such a normal form. 
   \item $i\Omega$ has a null eigenvector with nonzero eigenvalue:
    \begin{center} 
    \begin{tabular}{lll} 
     when $n$ is even, 
     & 
     $\tilde \Omega = a_1 S \oplus a_2 J \oplus \cdots \oplus a_m J \oplus 0$, 
     & 
     with $m = n / 2$; \\ 
     when $n$ is odd, 
     &
     $\tilde \Omega = a_1 S \oplus a_2 J \oplus \cdots \oplus a_m J$, 
     &
     with $m = (n+1)/2$.  
    \end{tabular}  
   \end{center} 
   Here, $a_1 > 0$ and $a_2 \geqslant \cdots \geqslant a_m \geqslant 0$.
   See \S 3.3.2 for an example. 
   \item $i\Omega$ has a null eigenvector with zero eigenvalue 
            but no timelike eigenvector:
    \begin{center} 
    \begin{tabular}{lll} 
     when $n$ is even, 
     & 
     $\tilde \Omega = a_1 T \oplus a_2 J \oplus \cdots \oplus a_m J$, 
     & 
     with $m = n / 2$; \\ 
     when $n$ is odd, 
     &
     $\tilde \Omega = a_1 T \oplus a_2 J \oplus \cdots \oplus a_m J \oplus 0$, 
     &
     with $m = (n-1)/2$.  
    \end{tabular}  
    \end{center} 
   Here, $a_1 > 0$ and $a_2 \geqslant \cdots \geqslant a_m \geqslant 0$.
   \S 3.3.3 exemplifies this normal form.
  \end{itemize}  
  In the above description, $J$, $S$ and $T$ denote the matrices
  \begin{displaymath}
   J = 
   \left(
    \begin{array}{cc}
      0  &  1  \\
     -1  &  0
    \end{array}
   \right) \, , 
   \quad \quad  
   S = 
   \left(
    \begin{array}{cc}
     0  &  1  \\
     1  &  0
    \end{array}
   \right) \, ,  
   \quad \quad 
   T = 
   \left(
    \begin{array}{ccc}
     0  &  1  &  0  \\
     1  &  0  &  1  \\
     0  & -1  &  0  \\
    \end{array}
   \right) \, . 
  \end{displaymath}  
  We declare this $\tilde \Omega$ to be the normal 
  form of $\Omega$.  It remains to determine how the strong convexity 
  criterion $| \tilde W | < 1$ constrains the parameters that describe these 
  normal forms.  This is where inequality $(2,-)_k$ of Theorem 7 comes into 
  play.  It reads: $(\tilde Qx + \tilde C) \cdot (\tilde Qx + \tilde C) - 
  (x \cdot \tilde C)^2 < |K| ( 1 - x \cdot x )$.

                     %.............................

 \subsubsection{When $i\Omega$ has a timelike eigenvector} 

  The type $(J)$ normal form $\tilde \Omega$ is derived in \S 10.3.4.  The 
  corresponding Killing field is given by $\tilde C = 0$ and 
  \begin{center}
   \begin{tabular}{lll}
    when $n$ is even,  & 
    $- \tilde Q = a_1 J \oplus \cdots \oplus a_m J$,  & 
    with $m=n/2$;   \\
    when $n$ is odd,  &
    $- \tilde Q = a_1 J \oplus \cdots \oplus a_m J \oplus 0$, &
    with $m = (n-1)/2$.  
   \end{tabular} 
  \end{center} 
  Here, $a_1 \geqslant a_2 \geqslant \cdots \geqslant a_m \geqslant 0$.
  Because $\tilde W(0) = 0$, the criterion $| \tilde W | < 1$ will always be 
  satisfied in some neighborhood of the origin.  Therefore the moduli space 
  is parametrized by 
  \begin{displaymath}
   a_1 \geqslant \cdots \geqslant a_m \geqslant 0 \, .
  \end{displaymath}

                   %...................................

 \subsubsection{When $i\Omega$ has a null eigenvector with nonzero eigenvalue} 

  The type $(S)$ normal form $\tilde \Omega$ is given in \S 10.3.5.  The 
  associated Killing field has data $\tilde C = (a_1, 0, \ldots, 0)$ and 
  \begin{center}
   \begin{tabular}{lll} 
    when $n$ is even,  & 
    $- \tilde Q = 0 \oplus a_2 J \oplus \cdots \oplus a_m J \oplus 0$, & 
    with $m=n/2$; \\
    when $n$ is odd,  &
    $- \tilde Q = 0 \oplus a_2 J \oplus \cdots \oplus a_m J$, &
    with $m=(n+1)/2$.  
   \end{tabular} 
  \end{center}
  Here, $a_1 > 0$ and $a_2 \geqslant \cdots \geqslant a_m \geqslant 0$.  
  The condition $| \tilde W | < 1$ is equivalent to 
  \begin{displaymath}
   a_1^2 ( 1 - x_1^2 ) \ + \ 
   \sum_{j=2}^{m} a_j^2 ( x_{2j-2}^2 + x_{2j-1}^2 ) 
   \ < \ |K| ( 1 - x \cdot x ) \, .
  \end{displaymath} 
  In particular, we must have 
  $a_1^2 ( 1 - x_1^2 ) < |K| ( 1 - x \cdot x )$.  This forces 
  $a_1 < \sqrt{ |K| }$.  Conversely, as long as $a_1$ satisfies this bound, 
  we shall have $| \tilde W | < 1$ on a neighborhood of the origin.  Hence 
  the moduli space is parametrized by 
  \begin{displaymath}
   \sqrt{|K|} > a_1 > 0 \ , \quad 
   a_2 \geqslant \cdots \geqslant a_m \geqslant 0 \, .
  \end{displaymath}

                   %..................................

  \subsubsection{When $i\Omega$ has a null eigenvector with zero eigenvalue 
                but no timelike eigenvector} 

  For this case, the normal form $\tilde \Omega$ is of type $(T)$ and is  
  determined in \S 10.3.6.  The corresponding Killing field $\tilde W$ is 
  specified by $\tilde C = (a_1, 0, \ldots, 0)$ and 
  \begin{center}
   \begin{tabular}{lll} 
    when $n$ is even,  & 
    $- \tilde Q = a_1 J \oplus \cdots \oplus a_m J$, & 
    with $m=n/2$;  \\
    when $n$ is odd,  &
    $- \tilde Q = a_1 J \oplus \cdots \oplus a_m J \oplus 0$, &
    with $m=(n-1)/2$.  
   \end{tabular} 
  \end{center} 
  Here, $a_1 > 0$ and $a_2 \geqslant \cdots \geqslant a_m \geqslant 0$.  
  Given this data, $|\tilde W | < 1$ precisely when 
  \begin{displaymath}
   a_1^2 ( 1 - x_2 )^2 \ + \
   \sum_{j = 2}^m  a_j^2 ( x_{2j-1}^2 + x_{2j}^2 ) \ < \
   |K| ( 1 - x \cdot x ) \, . 
  \end{displaymath}
  Consider a point $x$ of the type $(0, x_2, 0, \ldots, 0)$.  For this $x$, 
  the inequality takes the form $a_1^2 ( 1 - x_2 ) < |K| ( 1 + x_2 )$, which 
  always holds provided that $x_2$ is sufficiently close to 1.  Continuity 
  then extends the inequality to a neighbourhood of that $x$.  Thus strong 
  convexity does not impose any constraint on the $a_i$.  We conclude that 
  the moduli space is parametrized by 
  \begin{displaymath}
   a_1 > 0 \ , \quad a_2 \geqslant \cdots \geqslant a_m \geqslant 0 \, .
  \end{displaymath}

 		    %..............................

   \subsubsection{The moduli space for $K<0$}
   
   Unlike those of positive and zero flag curvature, Randers spaces of 
   negative constant flag curvature may arise in two different fashions, 
   corresponding to the cases $\sigma \not= 0$ and $\sigma = 0$.  
   Since $\sigma$ is invariant under isometries (\S 7.1), it makes sense to 
   talk about the local isometry classes, and hence the moduli spaces, for 
   these two families. 
   \begin{itemize} 
    \item 
     Zermelo navigation on Euclidean space under an infinitesimal homothety 
     with $\sigma \not= 0$ produces a metric with flag curvature 
     $K = -\onesixteenth \sigma^2$.  The local moduli space of these metrics 
     is parametrized in \S 7.3.2.  
    \item 
     For $\sigma = 0$, the perturbation of a Riemannian space form of negative 
     sectional curvature $K$ by an infinitesimal isometry generates a metric 
     with flag curvature $K$.  These spaces are parametrized, up to local 
     isometry, in \S 7.4.1--7.4.3.  
   \end{itemize} 
   Together the Euclidean and hyperbolic parametrizations provide a complete 
   description of the local isometry classes.
   \begin{proposition}
    The local isometry moduli space of $n$-dimensional strongly convex Randers 
    spaces of constant flag curvature $K<0$ is parametrized by 
    $a = (a_1, \ldots, a_m) \in \mathbb{R}^m$ as follows.
    \begin{itemize}
     \item[$(e)$] $h$ is the Euclidean metric and $\sigma = \pm 4 \sqrt{|K|}$.
                  The parameter space is         
                  \begin{displaymath}
                   a_1 \geqslant \cdots \geqslant a_m \geqslant 0 \, , 
                  \end{displaymath} 
                  where $m=n/2$ when $n$ is even, and $m=(n-1)/2$ when $n$ 
                  is odd.  These metrics may not be extended to all of 
                  $\mathbb{R}^n$.
     \item[$(k)$] $h$ is the Klein metric.  The parameter space is the 
                  disjoint union of three sets.
                  \begin{itemize} 
                   \item[$\circ$] When $n$ is even, $m=n/2$ and the three 
                                   sets are
                         \begin{eqnarray*}
	                  & a_1 \geqslant \cdots \geqslant a_m \geqslant 0 \, ,
                            & \\
	                  & \sqrt{|K|} > a_1 > 0 \, ,\ 
                            a_2 \geqslant \cdots a_m \geqslant 0\, , 
                            & \\
 	                  \hbox{and} 
                          & a_1 > 0 \, , \ 
                            a_2 \geqslant \cdots \geqslant a_m \geqslant 0 \, .
	                 \end{eqnarray*}
                   \item[$\circ$] When $n$ is odd, $m=(n+1)/2$ and the three 
                                  sets are
                         \begin{eqnarray*}
	                  & a_1 \geqslant \cdots \geqslant a_m = 0 \, ,  
                            & \\
	                  & \sqrt{|K|} > a_1 > 0 \, ,\ 
                            a_2 \geqslant \cdots a_m \geqslant 0\, , 
                            & \\
 	                  \hbox{and} 
                          & a_1 > 0 \, , \ 
                            a_2 \geqslant \cdots \geqslant a_m = 0 \, .
	                 \end{eqnarray*}
                  \end{itemize}
                  The globally defined metrics are parametrized by $a = 0$.
    \end{itemize}
   \end{proposition}

%------------------------------------------------------------------------------

\section{Discussion of projective flatness} 

 Let $M$ be an $n$-dimensional differentiable manifold. A metric on $M$ is 
 said to be projectively flat if $M$ can be covered by coordinate 
 charts in which the geodesics of the metric are straight lines.  For 
 Riemannian metrics, Beltrami's theorem says that the only projectively flat 
 ones are those with constant sectional curvature.  There are Finsler metrics 
 of constant flag curvature which are not projectively flat; see for 
 example \cite{Br96, Br02} and \cite{BS02}. Thus Beltrami's theorem does not 
 extend to the Finslerian setting.

                  %--------------------------------------

 \subsection{Douglas' theorem} 
 
  A theorem due to Douglas \cite{D28} states that a Finsler metric $F$ is 
  projectively flat if and only if two special curvature tensors are zero. 
  The first is the Douglas tensor. The second is the projective Weyl tensor 
  for $n \geqslant 3$, and the Berwald--Weyl tensor \cite{B47a} for $n=2$.  
  (The projective Weyl tensor automatically vanishes when $n=2$, thereby 
  predicating the need for a different invariant in that dimension.)  A 
  complete statement of Douglas' theorem can be found on p.144 of 
  \cite{Ru59}. 

  The projective Weyl tensor vanishes if and only if the flag curvatures of 
  $F$ have no dependence on the transverse edges (but can possibly depend on 
  the position $x$ and the flagpole $y$); see \cite{Sz77, M80}.  The 
  Berwald--Weyl tensor is defined for all $n$, though only relevant in 
  Douglas' theorem when $n=2$.  Shen has shown that, at least for the $n=2$ 
  case, this tensor is zero whenever the Ricci scalar and the S-curvature 
  (divided by $F$) are both constant.  (The Ricci scalar is the sum of 
  $n-1$ appropriately chosen flag curvatures.) 
    
                  %--------------------------------------

 \subsection{Specializing to Randers metrics} 

  For Randers metrics of constant flag curvature, we see that the projective 
  Weyl tensor vanishes.  Since the flag curvature is constant, so is the Ricci 
  scalar.  Moreover, for such metrics the S-curvature (divided by $F$) is     
  $\fourth (n+1) \sigma$, where $\sigma$ is the constant we encountered in 
  \S 4.1; see \cite{CS03}.  Thus the Berwald--Weyl tensor in two dimensions 
  is zero as well.  

  According to \cite{BaM97}, a Randers metric $F$ has vanishing Douglas 
  tensor if and only if the 1-form $b := b_i \, dx^i$ is closed.  
  Let $W^\flat$ denote the 1-form $W_i \, dx^i$, where $(h,W)$ is the Zermelo 
  navigation data of $F$. Using the equation $\mathcal L_W h = - \sigma \, h$ 
  with constant $\sigma$, it can be checked that the 2-forms 
  $\txtcurl := - d b$ (\S 4.1) and $\mathcal C := - d W^\flat$ (\S 4.3.2) 
  are related through $\txtcurl^{ij} = - \lambda \, \mathcal C^{ij}$, where 
  $\lambda := 1 - |W|^2$ is positive because of strong convexity (\S 2.2). 
  In particular, $db = 0 \Leftrightarrow dW^\flat = 0$, 
  whenever the above $\mathcal L_W$ equation holds. If the Randers 
  metric $F$ has constant flag curvature, then Theorem 3 (\S 4.4) avails us 
  of this $\mathcal L_W$ equation; in that case, the vanishing of the Douglas 
  tensor is equivalent to the condition $d W^\flat = 0$. 
 
                  %-------------------------------------

 \subsection{Projectively flat Randers metrics of constant flag curvature} 

  By virtue of Douglas' theorem, we see that a Randers metric $F$ of constant 
  flag curvature is projectively flat if and only if the 1-form 
  $W^\flat$ is closed, namely, $\partial_{x^j} W_i - \partial_{x^i} W_j = 0$. 
  \begin{itemize} 
   \item Suppose $F$ is obtained by perturbing the Euclidean metric. Using the 
         formula for $W_i$ given in the proof of Proposition 5, we see that 
         $W^\flat$ is closed if and only if $(Q_{ij})$ is the zero matrix. 
   \item Suppose $F$ is obtained by perturbing the standard sphere or the 
         Klein model. Using the formula for $W_i$ given in the proof of 
         Proposition 6, we find that $W^\flat$ is closed if and only if 
         $W$ is identically zero. The absence of non-trivial closed Killing 
         fields (namely parallel fields) on non-flat Riemannian space forms 
         is consistent with deRham's decomposition theorem.  
  \end{itemize} 
  
  The above information, together with our classification (Theorem 7), 
  tells us the following.  Up to local isometry, projectively flat strongly 
  convex {\it non}-Riemannian Randers metrics of constant flag curvature $K$ 
  comprise exactly two camps. 
  \begin{itemize} 
   \item[(1)] {\bf $K=0$:} 
              Zermelo navigation on Euclidean space with a constant vector 
              field $W = C$ satisfying $0< |C| < 1$.  These are the 
              (locally) Minkowski spaces; see \S 3.2.3.  A rotation can be 
	      used to transform $W$ into 
              $(0, \dots, 0, |C|)$ without causing the Minkowski metric in 
              question to leave its local isometry class.  Thus $|C|$ 
              parametrizes the 1-dimensional moduli space, which is the open 
              unit interval.  Excluding Randers metrics which are Riemannian 
              from Proposition 9 gives the same conclusion. 
   \item[(2)] {\bf $K<0$:} 
              Zermelo navigation on Euclidean space with 
              $W = - \half \sigma x + C$, $\sigma = \pm 4 \sqrt{|K|}$, 
              and $C \cdot C + \sigma x \cdot (\fourth \sigma x - C) < 1$. 
              This camp includes the Funk metric of \S 3.2.2.   
              A translation transforms $W$ into 
              $\tilde W = - \half \sigma \tilde x$.  By \S 7.3.2 and \S 7.1, 
              the corresponding Randers metrics $F$ and $\tilde F$ share the 
              same local isometry class.  Closer examination of $\tilde F$ 
              reveals that it is a $\tilde x$-scaled variant of the Funk 
              metric, one which lives on the open ball of radius 
              $1 / (2 \sqrt{|K|})$ centered at the origin of $R^n$.  In 
              particular, the moduli space consists of only one point, as 
              predicted by case $(e)$ of Proposition 10 (with all $a_i$ set 
              to zero because $Q=0$ here).  
  \end{itemize} 
  As a corollary of this itemization, every projectively flat strongly convex 
  Randers metric of constant positive flag curvature must be locally isometric 
  to a Riemannian standard sphere. 
  
  We see from the table in \S 5.4 that among the examples in \S 3, only 
  3.2.2 and 3.2.3 are projectively flat.   

                %---------------------------------------

 \subsection{Comments}  

  The above conclusions about projectively flat Randers metrics $F$ of 
  constant flag curvature is consistent with the main result of \cite{S02a}. 
  However, other than the fact that the two papers use totally different 
  methods, there are further distinctions.  Here, the $K<0$ camp has simple 
  navigation data $(h, W)$, where $h$ is the Kronecker delta; but the 
  resulting $F$, when generated with \S 2.1.3, shows a certain amount of 
  complexity.  In \cite{S02a}, a simple expression is derived for $F$ in the  
  $K<0$ camp; but, upon the use of \S 2.3 to recover the navigation 
  data $(h, W)$, we find that $h$, though isometric to the Euclidean metric, 
  takes on a somewhat complicated form.  Also, the fact that the moduli space 
  for $K<0$ consists of a single point is not manifest in \cite{S02a}. 
 
%------------------------------------------------------------------------------

\section{Restricting to the $\theta = 0$ family} 

 Recall the tensor $\theta_i := b^s \, \txtcurl_{si}$ encountered in \S 4.1.  
 Strongly convex Randers metrics of constant flag curvature and satisfying 
 the additional condition $\theta = 0$ have previously been characterized by 
 the {\it corrected} Yasuda--Shimada theorem.  See \cite{BR03, MS02} for 
 details and references therein. 

               %-----------------------------------------

 \subsection{Necessary and Sufficient conditions for $\theta=0$}
 
  It can be shown (using the machinery in \cite{BR04}) that the tensor 
  $\theta$ for Randers metrics of constant flag curvature has the navigation 
  description $(1-|W|^2) \theta_j = (|W|^2)_{:j} + \sigma W_j$.  Since our 
  Randers metrics are always presumed to be strongly convex ($|W|<1$), we see 
  that 
  \begin{displaymath} 
   \theta = 0 
   \quad 
   \Leftrightarrow  
   \quad 
   (|W|^2)_{:j} + \sigma W_j = 0 \, . 
  \end{displaymath}  

               %........................................

 \subsubsection{The Euclidean case}

  When $h$ is the flat Euclidean metric, $W = - \half \sigma x + Qx + C$ 
  according to Proposition 5.  The equation $(|W|^2)_{:j} + \sigma W_j = 0$ is 
  polynomial in the local coordinates $(x^i)$.  By considering the 
  coefficients of this polynomial, one can establish that $\theta = 0$ 
  if and only if 
  \begin{itemize}
    \item $Q= 0$ when $\sigma \not= 0$;  
    \item $Q^2 = 0$ and $QC = 0$ when $\sigma = 0$.
  \end{itemize}
  It is clear, from the normal form $\tilde Q$ (\S 10.2) of $Q$, that 
  $Q^2 = 0$ if and only if $Q = 0$.  Hence the two cases can be unified into a 
  single criterion $Q=0$, which is in turn equivalent to the 1-form 
  $W^\flat := W_i dx^i$ being closed (\S 8.3).  We conclude that, for strongly 
  convex constant flag curvature Randers metrics which are generated by 
  navigating on Euclidean $\mathbb R^n$ under the influence of an 
  infinitesimal homothety $W$, 
  \begin{center} 
   $\theta = 0$ \quad if and only if \quad $d W^\flat = 0$. 
  \end{center} 
  Such metrics are precisely the projectively flat ones enumerated in \S 8.3. 
  It is worth recollecting (\S 8.2) that in the present context, 
  $d W^\flat = 0$ is equivalent to $d b = 0$. 

              %.........................................
  
 \subsubsection{The spherical and Klein models}

  When $h$ is either the spherical or hyperbolic metric, 
  $\sigma$ must vanish (\S 4.4), and we see that 
  \begin{center}
   $\theta = 0 
    \quad \Leftrightarrow \quad 
    (|W|^2)_{:j} = 0
    \quad \Leftrightarrow \quad  
    |W| \ \textrm{is constant}$.
  \end{center}
  Proposition 6 says that 
  $W_i = (Q_{ij} x^j + C_i) / \{ |K| (1 + \psi x \cdot x) \}$, where 
  $\psi := K / |K|$.  Consequently the expression $(|W|^2)_{:j}$ 
  is rational in the local coordinates $(x^i)$.  Hence $\theta = 0$ if 
  and only if the polynomial numerator does, which ultimately leads to the 
  following necessary and sufficient conditions: 
  \begin{displaymath}       
   QC = 0 
   \quad \textrm{and} \quad 
   Q^2 = \psi \, (C C^t - |C|^2 I_n) \, .
  \end{displaymath}
  Here, $C$ is a column and $C^t$ is a row.  

  The above equations are invariant in form under any orthogonal 
  transformation $R \in O(n)$.  Indeed, multiplying each term by $R$ on the 
  left and $R^t$ on the right, those equations become $\tilde Q \tilde C = 0$ 
  and $\tilde Q^2 = \psi \, (\tilde C \tilde C^t - | \tilde C |^2 I_n)$,  
  where $\tilde Q = R Q R^t$, $\tilde C = RC$.  
  \begin{itemize} 
   \item Therefore, without any loss of generality, we may assume that $Q$ is 
         already in the normal form derived in \S 10.2.  Namely, 
         \begin{displaymath} 
          Q = q_1 J \oplus \cdots q_k J \oplus 0_{n-2k}, 
          \quad \textrm{with} \quad 
          q_1 \geqslant \cdots \geqslant q_k > 0 .
         \end{displaymath}   
   \item With this $Q$, the equation $QC=0$ can be solved immediately to find 
         that the first $2k$ components of $C$ are zero.  Its remaining 
         components can be transformed by any $r \in O(n-2k)$ without altering 
         $Q$.  Thus we may assume that the column vector $C$ which solves 
	 $QC=0$ has the simplified form 
         \begin{displaymath} 
          C = (0, \dots, 0, |C|) . 
         \end{displaymath} 
  \end{itemize}    
  
  We now substitute the displayed $Q$ and $C$ into the equation 
  $Q^2 = \psi \, (C C^t - |C|^2 I_n)$.  The outcome reads 
  \begin{displaymath}
   ({}^*) \quad \quad 
   -q_1^2 I_2 \oplus \cdots \oplus -q_k^2 I_2 \oplus 0_{n-2k} 
   = - \psi \, |C|^2 \, I_{n-1} \oplus 0 \, , 
  \end{displaymath}
  where $I_j$ denotes the $j \times j$ identity matrix. 
  \begin{itemize}
   \item By inspection, all the $q_i$ are zero if and only if $|C|=0$.  In 
         other words, $Q=0 \Leftrightarrow C=0$.  The Killing field  
         corresponding to $Q=0$, $C=0$ is $W=0$.  In that case the associated 
         Randers metric is simply the original Riemannian space form $h$.
   \item It remains to examine the scenario in which neither $Q$ nor $C$ is 
         identically zero.  Since all the $q_i$, as well as $|C|$, are 
         nonzero, equation $({}^*)$ forces three restrictions.   
         \begin{itemize}
          \item[(1)] $\psi := K / |K| = 1$, hence $K>0$ and $h$ must be 
                     the spherical metric.  
          \item[(2)] $q_1 = \cdots = q_k = |C|$. 
          \item[(3)] $2k = n-1$; equivalently, $n=2k+1$ is odd.
         \end{itemize}
         Up to local isometry, the strongly convex Randers metric in question 
         must have arisen from navigation on an {\it odd} dimensional sphere, 
         under the influence of a one parameter family of winds $W$. 
  \end{itemize} 
  We hasten to reiterate that these restrictions are obtained from local 
  considerations only; global conditions are not needed in their derivation. 

                     %..................................

 \subsection{Refined conclusions of the corrected Yasuda--Shimada theorem}

  Taken together, the previous subsections and \S 8.3 allow us to enumerate 
  all strongly convex Randers metrics with constant flag curvature $K$ and 
  $\theta = 0$.  They are obtained by Zermelo navigation on Riemannian space 
  forms $h$, subject to the influence of appropriate winds $W$ which satisfy 
  $|W| < 1$.  The {\it non}-Riemannian ones are as follows. 
  \begin{itemize}
   \item When $K < 0$: $h$ is the flat metric on Euclidean $\mathbb R^n$, 
                       and $W = - \half \sigma x + C$, with 
                       $\sigma = \pm 4 \sqrt{|K|}$.  As explained in 
                       \S 8.3, the resulting Randers metric is locally 
                       isometric to a position-scaled variant of the Funk 
                       metric, one which is generated by 
                       $\tilde W = - \half \sigma \tilde x$ and lives on the 
                       open ball of radius $1 / (2 \sqrt{|K|})$. 
   \item When $K = 0$: $h$ is the flat metric on Euclidean $\mathbb R^n$, 
                       and $W = C \not= 0$, living on the domain where 
                       $|C| < 1$.  We saw in \S 8.3 that up to local isometry, 
                       this family, which consists of locally Minkowski 
                       metrics, is parametrized by a single parameter $|C|$. 
   \item When $K > 0$: $h$ is $1 / K$ times the standard metric on the 
                       unit sphere $S^n$, with $n = 2k+1$ odd.  The wind $W$ 
                       is given in projective coordinates (\S 5.3, \S 6.2.1) 
                       as $Qx + C + (x \cdot C)x$, where $Q$ and $C$ are 
                       specially related on account of $\theta = 0$. 
                       In fact (\S 9.1.2), there is an $R \in O(n)$ such that 
                       $\tilde C := RC = (0, \dots, 0, |C|)$ and 
                       $\tilde Q := RQR^t = 
                        |C| (J \oplus \cdots \oplus J) \oplus 0$, 
                       respectively.   This is equivalent to conjugating the 
                       matrix representative (\S 7.2) of $W$ by the element 
                       $1 \oplus R$ in the isometry group of $h$.  Thus 
                       (\S 7.1) the Randers metric generated by 
                       $\tilde W := 
                        \tilde Qx + \tilde C + (x \cdot \tilde C)x$ lies in 
                       the same isometry class as that from $W$.  
                       Applying the analysis in \S 7.2.2 to 
		       $\tilde W$, we see that strong convexity mandates 
		       $|C| < \sqrt{K}$, which as a bonus (\S 7.2.3) ensures 
		       that the metric is global on $S^n$.  Thus, up to 
		       isometry, there is only a one parameter family (indexed 
		       by $|C|$) of {\it non}-Riemannian strongly convex 
                       Randers metrics with constant flag curvature $K$ and 
                       $\theta = 0$ on the odd dimensional spheres.  By 
                       contrast, no such metric exists on the even dimensional 
                       spheres, regardless of whether it is locally or 
                       globally defined.   
  \end{itemize}    

  Strongly convex non-Riemannian Randers metrics with constant flag curvature 
  $K$ and $\theta = 0$ are characterized by the corrected Yasuda--Shimada 
  theorem \cite{BR03, MS02}.  The conclusion for the $K=0$ case is as 
  described above.  For nonzero $K$, the characterization is in terms of 
  coupled systems of nonlinear partial differential equations.  Our discussion 
  above may be viewed as a complete list of solutions to those partial 
  differential equations. 
    
  Bejancu--Farran \cite{BF02,BF03}, with the help of the corrected 
  Yasuda--Shimada theorem, have recently established a bijection between 
  Sasakian space forms of constant $\phi$-sectional curvature $c \in (-3,1)$, 
  and Randers metrics of constant flag curvature $K=1$ with $\theta = 0$.  In 
  the course of their study they showed that the underlying manifold $M$ must 
  be of odd dimension, and is necessarily diffeomorphic to a sphere when it is 
  simply connected and complete with respect to $a$.  These results can be
  made equivalent to what we have described for the $K > 0$ case.  Our 
  $\theta$ is denoted by $\beta$ in the Bejancu--Farran papers, and their $c$ 
  is $1 - 4 \| b \|^2$ in our notation.  
                  
  It is worth mentioning here that all spheres, of both odd and {\it even} 
  dimensions, admit a wealth of non-Riemannian globally defined Randers metrics
  of constant positive flag curvature, provided that the restriction 
  $\theta = 0$ is lifted.  
  
  Here is a straightforward example on $S^4$. 
  Following the treatment of \S 6.2 we let $p = (p^0, p^1, p^2, p^3, p^4)$ 
  denote the canonical coordinates on $\mathbb{R}^{5}$.  The infinitesimal 
  rotation 
  \begin{displaymath}   
   W(p) = \tau ( - p^2 \partial_{p^1} + p^1 \partial_{p^2} ) \, , 
   \quad \tau \ \textrm{constant}, 
  \end{displaymath} 
  restricts to a globally defined Killing field on the standard unit sphere 
  $S^4$.  As long as $|\tau| < 1$, we have $|W| < 1$ on the entire sphere.  
  Hence $W$ induces a globally defined, strongly convex Randers metric with 
  constant flag curvature $+1$ on $S^4$.

  Notice, however, that $\theta \not= 0$.  This is immediate from the equation 
  displayed at the beginning of \S 9.1.2.  It says that $\theta$ vanishes if 
  and only if $|W|$ is constant.   The norm of our $W$ is certainly not 
  constant.  Hence $\theta$ is nonzero.

%------------------------------------------------------------------------------

\section{Appendix: Some Lie theory}
 
 Recall from \S 7.1 that the symmetry/isometry groups $G$ (of the Riemannian 
 space forms) act on the Lie algebras of infinitesimal homotheties, via the 
 adjoint action $Ad$.  Our analysis of the moduli space (\S 7) of constant 
 flag curvature Randers metrics requires detailed knowledge of each $Ad$ 
 orbit, in order to pinpoint a specific representative in the fundamental 
 Weyl chamber.  

 Though the Lie theory for the orthogonal group is well known, it is invoked 
 so many times in the paper, and in several different contexts, that we feel 
 obligated to at least set the notation and state the facts (\S 10.2).  In the 
 non-compact case $G = O_+(1,n)$, the orthochronous Lorentz group, the 
 information we need is less standard, and is typically not in a form that we 
 could use without substantial modification or synthesis.  Since this 
 information plays such a pivotal role in our geometrical conclusions, we are 
 compelled to sketch a cohesive account (\S 10.3), which takes up the bulk of 
 the Appendix. 

 Finally, we have chosen to present this material in matrix language for the 
 sake of concreteness. 

                %.........................................

 \subsection{Scalar products and the ``perp argument"} 
  
  By a {\it scalar product} on any complex vector space $\mathcal V$, we mean 
  a pairing $\langle \, , \, \rangle$ which is $\mathbb C$-linear in the first 
  factor, satisfies $\langle u, v \rangle = \overline{\langle v, u \rangle}$, 
  and is non-degenerate (namely, if $\langle u, v \rangle = 0$ for all $v \in 
  \mathcal V$, then $u$ must vanish).  Inner products are simply positive 
  definite scalar products.  For example, if $E$ is the diagonal matrix 
  $-1 \oplus I_n$, then $\langle u, v \rangle := u^t E \overline v$ is a 
  scalar product on $\mathbb C^{1+n}$, whereas replacing that $-1$ by $+1$ 
  gives the canonical inner product $u^t \overline v$ on $\mathbb C^{1+n}$. 
  A vector $v$ is said to be spacelike, null, or timelike, respectively, 
  if $\langle v, v \rangle$ is positive, zero, or negative. 

  Let $\mathcal W$ be any subspace of a scalar product space $\mathcal V$.  
  Its perp $\mathcal W^\perp$ is $\{ v \in \mathcal V : 
  \langle v, w \rangle = 0 \ \textrm{for all} \ w \in \mathcal W \}$.  The 
  restriction of $\langle \, , \, \rangle$ to $\mathcal W^\perp$ may  
  fail to be non-degenerate when $\mathcal W$ contains a null vector. 
  For instance, in $\mathbb C^{1+2}$ with $E = \textrm{diag}(-1,1,1)$, if 
  $\mathcal W = \textrm{span} \{ (1,1,0) \}$, then $\langle \, , \, \rangle$ 
  is degenerate on $\mathcal W^\perp = \textrm{span} \{ (1,1,0), (0,0,1) \}$.  
  On the other hand, if $\mathcal W = \textrm{span} \{ (1,1,0), (1,-1,0) \}$, 
  then non-degeneracy holds on 
  $\mathcal W^\perp = \textrm{span} \{ (0,0,1) \}$.  These examples illustrate 
  the following useful fact implicit in \cite{ON83}: 
  $\mathcal W$ admits a $\langle \, , \, \rangle$ orthonormal basis 
  $\Leftrightarrow$ 
  $W \cap W^\perp = \{ 0 \}$ 
  $\Leftrightarrow$ 
  the restriction of $\langle \, , \, \rangle$ to $\mathcal W^\perp$ is 
  non-degenerate, in which case it defines a scalar product there. 

  Let $A$ be a self-adjoint linear operator on the scalar product space 
  $\mathcal V$.  Suppose the subspace $\mathcal W$ is invariant under $A$.  
  Then so is $\mathcal W^\perp$, because 
  $\langle Av, w \rangle = \langle v, Aw \rangle$.  Hence the restriction of 
  $A$ to $\mathcal W^\perp$ makes sense.  If, in addition, 
  $\langle \, , \, \rangle$ is non-degenerate on $\mathcal W^\perp$, then the 
  restricted $A$ is again operating on a scalar product space, albeit a 
  smaller one.  We shall repeatedly invoke this ``perp argument". 
 
                %..........................................

 \subsection{A compact case: normal form for skew-symmetric real matrices}

  Let $\Omega$ be any real $\ell \times \ell$ skew-symme\-tric matrix.  Then 
  $A := i \Omega$ is a self-adjoint linear operator on the inner product 
  space $\mathbb C^\ell$, with $\langle u, v \rangle := u^t \overline v$. 
  \begin{itemize} 
   \item For eigenvectors $z_1$, $z_2$ with eigenvalues $\lambda_1$, 
         $\lambda_2$, respectively, self-adjointness leads 
         to $\lambda_1 \langle z_1, z_2 \rangle  
             = \overline{\lambda_2} \, \langle z_1, z_2 \rangle$.  Thus  
         eigenvectors corresponding to distinct eigenvalues are 
         $\langle \, , \, \rangle$ orthogonal, and all eigenvalues of $A$ 
         are real.   
   \item Since $A = i \Omega$ where $\Omega$ is real, we have 
         $Az = \lambda z$ if and only if 
         $A \overline z = - \overline \lambda \, \overline z$.  Hence the 
         nonzero eigenvalues of $A$ must occur in pairs $\pm a$ ($a>0$), with 
         $\langle \, , \, \rangle$ orthogonal eigenvectors 
         $z$ and $\overline z$. 
   \item For each nonzero pair $\pm a$, the real vectors 
         $v := (z + \overline z)/2$ and $u := (z - \overline z)/(2i)$  
         satisfy $A u = - ia v$ and $A v = ia u$.  The 
         $\langle \, , \, \rangle$ orthogonality between $z$ and $\overline z$ 
         gives $\langle u, u \rangle = \langle v, v \rangle$ and 
         $\langle u, v \rangle = 0$.  Thus, the normalized versions 
         $\hat u$, $\hat v$ are orthogonal real unit vectors that still 
         satisfy $A \hat u = - ia \hat v$ and $A \hat v = ia \hat u$.
  \end{itemize} 
 
  Enumerate the nonzero eigenvalues of $A$, counted with multiplicity, as 
  $\pm a_1, \dots, \pm a_k$, where $a_1 \geqslant \dots \geqslant a_k > 0$.  
  Associated to $\pm a_1$ is the subspace $\mathcal W_1 := \textrm{span} 
  \{ \hat u_1, \hat v_1 \}$ which is invariant under the self-adjoint $A$.  
  Since $\hat u_1, \hat v_1$ are orthonormal, the perp argument (\S 10.1) 
  says that $A$ restricts to a self-adjoint linear operator on the inner 
  product space $\mathcal W_1^\perp$, and its largest eigenvalue pair on 
  $\mathcal W_1^\perp$ is $\pm a_2$.  Repeating this perp argument $k$ times, 
  we obtain a special orthonormal set of real vectors 
  $\{ \hat u_1, \hat v_1, \dots, \hat u_k, \hat v_k \}$ which are 
  $\langle \, , \, \rangle$ orthogonal to the $(\ell - 2k)$-dimensional 
  generalized null space of $A$.  

  If $2k < \ell$, then $A$ admits $0$ as an eigenvalue of algebraic 
  multiplicity $\ell - 2k$, with at least one eigenvector $\xi_{2k+1}$, say of 
  unit length.  This $\xi_{2k+1}$ also belongs to the null space of $\Omega$, 
  and hence can be chosen to be real.  Successive applications of the perp 
  argument generates an orthonormal set of real eigenvectors 
  $\{ \xi_{2k+1}, \dots, \xi_\ell \}$ for the eigenvalue $0$.  

  Consequently, for every eigenvalue of $A$, the algebraic and geometric 
  multiplicities agree.  With respect to the the real orthonormal basis 
  $\{ \hat u_1, \hat v_1, \dots, \hat u_k, \hat v_k, 
      \xi_{2k+1}, \dots, \xi_\ell \}$, the 
  matrix representation of $A$ is 
  $ia_1 J \oplus \cdots \oplus ia_k J \oplus 0_{\ell-2k}$, where  
  \begin{displaymath}
   J = \left(
   \begin{array}{rr}
   0  & 1 \\
   -1 & 0 \\ 
   \end{array}
   \right) \, .
  \end{displaymath}       
  Correspondingly, that of $\Omega$ is 
  $\tilde \Omega := a_1 J \oplus \cdots \oplus a_k J \oplus 0_{\ell-2k}$, 
  with $2k$ being its rank.  Suppressing the rank of $\Omega$, we see that 
  \begin{center} 
   \begin{tabular}{lll} 
    when $\ell$ is even, 
    & 
    $\tilde \Omega = a_1 J \oplus \cdots \oplus a_m J$, 
    & 
    with $m = \ell / 2$; \\ 
    when $\ell$ is odd,  
    &
    $\tilde \Omega = a_1 J \oplus \cdots \oplus a_m J \oplus 0$, 
    &
    with $m = (\ell-1)/2$; 
   \end{tabular}  
  \end{center} 
  where $a_1 \geqslant a_2 \geqslant \cdots \geqslant a_m \geqslant 0$. 
  This is the desired normal form of $\Omega$.  Note that 
  $\tilde \Omega = B^{-1} \Omega B$, where $B$ is the orthogonal matrix 
  whose columns are given by the vectors in our real orthonormal basis. 

  In terms of Lie theory, a skew-symmetric matrix $\Omega$ is an element in 
  the Lie algebra $\mathfrak o(\ell)$ of the orthogonal group $O(\ell)$.  The 
  fact that $\exp(a_i J)$ is the $2 \times 2$ rotation matrix with angle 
  $a_i$ tells us that $\exp(\tilde \Omega)$ lands in a maximal torus of 
  $O(\ell)$, and $\tilde \Omega$ itself belongs to a Cartan subalgebra 
  $\mathcal H$ of $\mathfrak o(\ell)$.  The condition 
  $a_1 \geqslant \cdots \geqslant a_m \geqslant 0$ singles out the fundamental 
  closed Weyl chamber of $\mathcal H$.  Setting $g = B^{-1}$, our arguments 
  show that every real skew-symmetric $\Omega$ can be $O(\ell)$-conjugated 
  into this closed Weyl chamber.

           %-----------------------------------------------------

 \subsection{A non-compact case} 

  We now use the tools and notation set up in \S 10.2 to derive a normal form 
  for elements of the Lie algebra $\mathfrak o(1,n)$. 

  Let $E$ denote the diagonal matrix $-1 \oplus I_n$.  The elements of 
  $\mathfrak o(1,n)$ are real $(n+1) \times (n+1)$ matrices $\Omega$ which 
  satisfy the condition $\Omega^t = - E \Omega E$; equivalently, $\Omega$ 
  has the defining form  
  \begin{displaymath}
   \Omega = 
   \left(
    \begin{array}{cc}
     0  &  C^t \\
     C  &  -Q  
    \end{array}
   \right)  \, ,  
  \end{displaymath} 
  where $Q$, $C$ are real, and $Q$ is $n \times n$ skew-symmetric.  The 
  Lorentz group $O(1,n)$ is a non-compact Lie group of which 
  $\mathfrak o(1,n)$ is the Lie algebra.  Elements of $O(1,n)$ are real 
  $(n+1) \times (n+1)$ matrices $g$ such that $g^{-1} = E g^t E$.  Thus the  
  first column of $g$ is a time-like unit vector and the remaining columns are 
  space-like unit vectors.  In particular, the top left entry of $g$ satisfies 
  $(g^0{}_0)^2 \geqslant 1$.  For reasons that will be made clear later, our 
  interest is in the orthochronous subgroup $G := O_+(1,n)$, for which 
  $g^0{}_0 \geqslant 1$.  

                      %.............................

 \subsubsection{An available simplification} 

  Our goal here is to select a simplest representative along the $G$ adjoint 
  orbit of $\Omega$.  To that end, we first invoke \S 10.2 to find an  
  element $R \in O(n)$ such that 
  $RQR^{-1} = q_1 J \oplus \cdots \oplus q_h J \oplus 0_{n-2h}$, where 
  $q_1 \geqslant \cdots \geqslant q_h > 0$.  This has the effect of 
  changing $C$ to $RC$.  Next, we use an element $r \in O(n-2h)$ to 
  transform the last $n-2h$ components of $RC$ into $(0, \dots, 0, \xi)$ 
  without affecting its first $2h$ components.  In terms of matrix 
  conjugation, set $g_1 := 1 \oplus R$ and 
  $g_2 := 1 \oplus I_{2h} \oplus r$, then  
  $(g_2 g_1) \Omega (g_2 g_1)^{-1}$ has the simplified form 
  \begin{center}
   \begin{tabular}{l} 
   $
   \left(
    \begin{array}{cccc}
     0   &  D^t                                & 0   & \xi  \\
     D   & -(q_1 J \oplus \cdots \oplus q_h J) & 0   & 0    \\ 
     0   &  0                                  & 0   & 0    \\ 
     \xi &  0                                  & 0   & 0  
    \end{array}
   \right) \, . 
   $
   \end{tabular}  
  \end{center}  
  Here, $D$ is a column of $2h$ entries listed pairwise; namely,   
  $D = (D_1, \dots, D_h)$, with $D_j := [(RC)_j, (RC)_{j+1}]$. 
  Since $g_2 g_1 \in G$, the above matrix lies on the same $Ad$ orbit as 
  $\Omega$.  When necessary, we can use this simplified form 
  for $\Omega$ with no loss of generality. 

                     %...............................

 \subsubsection{Preliminaries about eigenvalues and eigenvectors} 

  Given any $\Omega \in \mathfrak o(1,n)$, the matrix $A := i \Omega$ is a 
  self-adjoint linear operator on the scalar product space $\mathbb C^{1+n}$, 
  with $\langle U, V \rangle := U^t E \overline V$.  Let $V = (v_0, v)$ be 
  any (possibly complex) eigenvector of $A$ with eigenvalue $\lambda$.  
  Direct verification tells us that: 
  \begin{itemize} 
   \item[(1)] both $AV = \lambda V$ and 
              $A \overline V = - \overline \lambda \ \overline V$ hold;  
   \item[(2)] we have $\lambda v_0 = i C^t v$ and 
              $\lambda v = i v_0 C - i Qv$;  
   \item[(3)] either $\lambda = 0$ or $v_0^2 = v^t v$.  This is equivalent to 
              $\lambda (v_0^2 - v^t v) = 0$, which follows from (2) and the 
              skew-symmetry of $Q$.  
  \end{itemize}  
  Suppose $\lambda \not= 0$, so that $v_0^2 = v^t v$ by (3) above. 
  \begin{itemize} 
   \item[(4)] Then $V$ is either space-like or null.  (Equivalently, all 
              time-like eigenvectors must have zero eigenvalue.)  This comes 
              about because $\langle V, V \rangle = - |v_0|^2 + |v|^2$  
              and $|v_0|^2 = |v^t v| = |(v, \overline v)| 
                   \leqslant |v| \, |\overline v| = |v|^2$, where the  
              Cauchy--Schwarz inequality is being applied to the canonical 
              inner product $( \, , \, )$ on $\mathbb C^n$. 
   \item[(5)] The space-like eigenvectors have real eigenvalues, which must 
              occur in pairs $\pm a$ $(a>0)$, with corresponding 
              $\langle \, , \, \rangle$ orthogonal eigenvectors 
              $V$, $\overline V$. 
              Note that $A$ being self-adjoint implies that 
              $\lambda \langle V, V \rangle 
               = \overline \lambda \ \langle V, V \rangle$, hence $\lambda$ is 
              real whenever $V$ is not null.  The rest follows from 
              $\lambda = a > 0$, item (1), and 
              $a \langle V, \overline V \rangle 
               = - a \langle V, \overline V \rangle$.  
   \item[(6)] The null eigenvectors have pure imaginary eigenvalues, and can 
              always be standardized into the form $V = (1,v)$ with $v$ 
              real.  Indeed, $V = (v_0, \tilde v)$ being nonzero and null 
              means that $|v_0|^2 = |\tilde v|^2$ with $v_0 \not= 0$; dividing 
              by $v_0$ gives $(1,v)$, where $v^t \overline v = |v|^2 = 1$. 
              Yet, (3) says that $v^t v = 1$.  Substituting 
              $v = \textrm{Re} \, v + i \textrm{Im} \, v$ 
              into these two equations gives $\textrm{Im} \, v = 0$.  Then (1) 
              tells us that $\lambda = - \overline \lambda$. 
  \end{itemize} 

                   %...................................

 \subsubsection{Categorizing the normal forms of $A = i\Omega$} 

  We first note that 
  \begin{center} 
   ``if $A$ has no timelike eigenvector,  
   then it must admit a null eigenvector." 
  \end{center} 
  Given the absence of timelike eigenvectors, suppose there were no null 
  eigenvectors either.  Then all eigenvectors of $A$ would have to be 
  spacelike.  Applying the perp argument (\S 10.1) $n$ times would produce a 
  $\langle \, , \, \rangle$ orthonormal basis $B$ which is entirely spacelike 
  (and which diagonalizes $A$).  With respect to $B$, the matrix of 
  $\langle \, , \, \rangle$ would be $I_{n+1}$ instead of $E = -1 \oplus I_n$, 
  contradicting the invariance of the index of $\langle \, , \, \rangle$.   

  Thus it is reasonable to split our derivation of the normal forms of $A$ 
  into three camps. 
  \begin{itemize} 
   \item When $A$ has a timelike eigenvector, the normal form is of 
         type $(J)$. 
   \item In the absence of timelike eigenvectors:   
         \begin{itemize}  
          \item[*] $A$ has a null eigenvector with nonzero eigenvalue, in 
                   which case its normal form is of type $(S)$. 
          \item[*] $A$ has a null eigenvector with eigenvalue zero; then its 
                   normal form is said to be of type $(T)$. 
         \end{itemize}  
  \end{itemize} 
  These are discussed separately in \S 10.3.4, \S 10.3.5 and \S 10.3.6.  After 
  those discussions, the following will be apparent: 
  (a) The three types of normal forms are mutually exclusive.  
  (b) The absence of timelike eigenvectors is essential for the type $(T)$ 
      normal form to surface.  
  (c) Having a null eigenvector with {\it nonzero} eigenvalue automatically 
      rules out timelike eigenvectors; hence the assumption about timelike 
      eigenvectors being absent is not needed in the type $(S)$ case.  

                   %...................................

 \subsubsection{In the presence of a timelike eigenvector for $A$} 

  Call this eigenvector $U$; by item (4) of \S 10.3.2, its eigenvalue must be 
  $0$.  This puts $U$ in the null space of $A$ and hence that of $\Omega$.  
  Since the latter is real, $U$ can be chosen real.  Being timelike, the first 
  component $u_0$ of $U$ cannot vanish.  Replace $U$ by $-U$ if necessary to 
  effect $u_0 > 0$, and scale $U$ to unit length. 
 
  Set $\mathcal U := \textrm{span} \{ U \}$.  Since $U$ is timelike, the 
  restriction of $\langle \, , \, \rangle$ to $\mathcal U^\perp$ is positive 
  definite.  Hence the analysis of $A_{|\mathcal U^\perp}$ reduces to the 
  compact case considered in \S 10.2.  So there is a real orthonormal 
  basis $B$ for $\mathcal U^\perp$, with respect to which 
  $A_{|\mathcal U^\perp}$ has the normal form 
  $i a_1 J \oplus \cdots \oplus i a_k J \oplus 0_{n-2k}$.  

  The collection $\mathcal B := \{ U \} \cup B$ is a real 
  $\langle \, , \, \rangle$ orthonormal basis which puts $\Omega$ into the 
  normal form $\tilde \Omega := 
  0 \oplus a_1 J \oplus \cdots \oplus a_k J \oplus 0_{n-2k}$, 
  with $a_1 \geqslant \cdots \geqslant a_k > 0$.  Since $u_0 > 0$, the 
  corresponding matrix $g := {\mathcal B}^{-1}$ belongs to $O_+(1,n)$. 
  Suppressing the rank of $\Omega$ gives the following:   
  \begin{center} 
   \begin{tabular}{lll} 
    when $n$ is even, 
    & 
    $\tilde \Omega = 0 \oplus a_1 J \oplus \cdots \oplus a_m J$, 
    & 
    with $m = n / 2$; \\ 
    when $n$ is odd, 
    &
    $\tilde \Omega = 0 \oplus a_1 J \oplus \cdots \oplus a_m J \oplus 0$, 
    &
    with $m = (n-1)/2$.  
   \end{tabular}  
  \end{center} 
  Here, $a_1 \geqslant a_2 \geqslant \cdots \geqslant a_m \geqslant 0$.
  
                   %...................................

 \subsubsection{When $A$ has a null eigenvector with nonzero eigenvalue} 

  Take any such null eigenvector and call it $X$.  According to item (6) of 
  \S 10.3.2, the eigenvalue in question has the form $ia$ with 
  $0 \not= a \in \mathbb R$, and $X$ can be chosen as $(1,x)$, where $x$ 
  is real and $|x|^2 = 1$.  Incidentally, item (2) of \S 10.3.2 
  characterizes $x$ by the equations $a = C^t x$ and $ax = C - Qx$.   

  There is in fact a companion real null eigenvector $Y$ with the standardized 
  form $(1,y)$, and which has eigenvalue $-ia$.  To see this, it suffices to 
  solve $-a = C^t y$ and $-ay = C - Qy$ for a real $y$.  The condition 
  $|y|^2 = y^t y = 1$ then follows from these two equations and the fact that 
  $a$ is nonzero.   
  
  Since $Q^t = -Q$, we can rewrite the second equation as 
  $y^t (Q+aI) = - C^t$.  Also, $Q+aI$ is invertible because the spectrum of 
  $Q$ is pure imaginary (\S 10.2).  Thus $y^t = - C^t (Q+aI)^{-1}$, which is 
  real because $Q$ and $C$ are.  Finally, with the help of the hypothesized 
  $x$, we have $C^t y = y^t C = y^t (Q+aI)x = - C^t x = -a$.  This proves that 
  the asserted $Y$ exists.  Further analysis, based on the self-adjointness of 
  $A$, shows that any standardized null eigenvector with nonzero eigenvalue 
  must be either $X$ or $Y$. 

  By interchanging $X$ with $Y$ if necessary, we may assume that $a>0$. 
  For later purposes, relabel it as $a_1$.  Define $U:=X+Y=(2,x+y)^t$ and 
  $V:=X-Y=(0,x-y)^t$.  Observe that: 
  \begin{itemize} 
   \item[*] $\langle U, U \rangle = 2(-1+x \cdot y) < 0$ and 
            $\langle V, V \rangle = 2(1-x \cdot y) > 0$;  
   \item[*] $U$ and $V$ are $\langle \, , \, \rangle$ orthogonal;  
   \item[*] $AU = ia_1 V$ and $AV = ia_1 U$.  Since 
            $| \langle U, U \rangle | = \langle V, V \rangle$, that pair of 
            equations remains valid for the normalized vectors $\hat U$ and 
            $\hat V$. 
  \end{itemize} 

  Set $\mathcal W := \textrm{span} \{ \hat U, \hat V \}$.  Since $\hat U$ is 
  timelike, the scalar product $\langle \, , \, \rangle$ becomes positive 
  definite on the $(n-1)$-dimensional $\mathcal W^\perp$, which  
  is invariant under the self-adjoint $A$.  In view of \S 10.2, 
  there is a real orthonormal basis $B$ for $\mathcal W^\perp$, with respect 
  to which the restricted $A$ has the normal form 
  $i a_2 J \oplus \cdots \oplus i a_k J \oplus 0_{n-1-2(k-1)}$.  

  The collection $\mathcal B := \{ \hat U, \hat V \} \cup B$ is a real 
  $\langle \, , \, \rangle$ orthonormal basis which puts $\Omega$ into the 
  normal form $\tilde \Omega := 
  a_1 S \oplus a_2 J \oplus \cdots \oplus a_k J \oplus 0_{n+1-2k}$, where  
  \begin{displaymath}
   S = 
   \left(
    \begin{array}{cc}
     0  &  1  \\
     1  &  0
    \end{array}
   \right) 
  \end{displaymath}   
  and $a_1 > 0$, $a_2 \geqslant \cdots \geqslant a_k > 0$.  Since the 
  first component of $\hat U$ is positive, the corresponding 
  matrix $g := {\mathcal B}^{-1}$ belongs to $O_+(1,n)$.  Suppressing the 
  rank of $\Omega$ gives the following:   
  \begin{center} 
   \begin{tabular}{lll} 
    when $n$ is even, 
    & 
    $\tilde \Omega = a_1 S \oplus a_2 J \oplus \cdots \oplus a_m J \oplus 0$, 
    & 
    with $m = n / 2$; \\ 
    when $n$ is odd, 
    &
    $\tilde \Omega = a_1 S \oplus a_2 J \oplus \cdots \oplus a_m J$, 
    &
    with $m = (n+1)/2$.  
   \end{tabular}  
  \end{center} 
  Here, $a_1 > 0$ and $a_2 \geqslant \cdots \geqslant a_m \geqslant 0$.

                   %..................................

 \subsubsection{When $A$ has a null eigenvector with zero eigenvalue 
                but no timelike eigenvector} 

  Items (6) and (2) of \S 10.3.2 tell us that such a null eigenvector $V$ can 
  always be standardized into the form $(1,v)$, where $v$ is real, 
  $v \cdot v = 1$, $Qv = C$, and $C \cdot v = 0$.  \S 10.3.1 says 
  there is no loss of generality in assuming that $Q$ and $C$ have already 
  been simplified to $q_1 J \oplus \cdots \oplus q_h J \oplus 0_{n-2h}$ 
  and $(D_1, \dots, D_h, 0, \dots, 0, \xi)$, respectively.  Here, 
  $q_1 \geqslant \cdots \geqslant q_h > 0$ and $D_j = [C_j, C_{j+1}]$.  The 
  hypothesized existence of $V$ implies that $Qv = C$ admits a solution.  
  Hence $C$ is in the range of $Q$ and $\xi$ must vanish.  The use of 
  $J^2 = - I$ solves the equation $Qv=C$ to give 
  \begin{displaymath}
   v = \left( \frac{-J D_1}{q_1}, \dots, \frac{-J D_h}{q_h}, 
              v_{2h+1}, \dots, v_n \right) \, .   
  \end{displaymath} 
  This $v$ automatically satisfies $C \cdot v = 0$ because of the 
  skew-symmetry of $J$, and its last $n-2h$ components are constrained by the 
  requirement $v \cdot v = 1$.  

  For further discussions, set  
  \begin{displaymath} 
   z := \left( \frac{-J D_1}{q_1}, \dots, \frac{-J D_h}{q_h}, 
               0, \dots, 0 \right) \, .
  \end{displaymath} 
  The null space $\mathcal N_1$ of $A = i \Omega$ consists of 
  eigenvectors $U = (u_0, u)$ with eigenvalue $0$, which are characterized 
  by $Qu = u_0 C$ and $C^t u = 0$.  Since $\Omega$ is real, $U$ may 
  be chosen to be real.  A calculation like the one above tells us that 
  $\mathcal N_1$ admits a basis 
  $\{ (1,z), (0,e_j), j = 2h+1, \dots, n \}$, 
  where $e_j$ has a $1$ in the $j$th entry, and $0$ elsewhere.  In particular, 
  $(1,z)$ is an eigenvector of $A$ with eigenvalue $0$.  

  If $(1,z)$ were not null, then the components $v_{2h+1}, \dots, v_n$ of 
  the hypothesized null eigenvector $(1,v)$ could not all be zero, whence  
  $|z|^2 < |v|^2 = 1$.  This would force the eigenvector $(1,z)$ to be 
  timelike, a scenario forbidden by our hypothesis.  Thus $(1,z)$ has to be 
  null.  

  Since $|JD_i| = |D_i|$, the condition $z \cdot z = 1$ is equivalent to 
  \begin{displaymath} 
   ({}^*) \quad\quad
   \frac{|D_1|^2}{q_1^2} + \cdots + \frac{|D_h|^2}{q_h^2} = 1 \, .  
  \end{displaymath} 
  Introduce the column vectors 
  \begin{displaymath} 
   z_1 := \left( \frac{D_1}{q_1^2}, \dots, \frac{D_h}{q_h^2}, 
                 0, \dots, 0 \right) , 
   \quad\quad  
   z_2 := \left( \frac{J D_1}{q_1^3}, \dots, \frac{J D_h}{q_h^3}, 
                 0, \dots, 0 \right) .
  \end{displaymath} 
  Let $\mathcal N_i$ be the null space of $A^i$ and abbreviate the vectors 
  $(1,z), (0,e_j), j = 2h+1, \dots, n$ collectively as $B_0$.  Then 
  \begin{center} 
   \begin{tabular}{l} 
    $\mathcal N_1 = \textrm{span} \{ B_0 \}$, 
    \quad 
    $\mathcal N_2 = \textrm{span} \{ (0,z_1), B_0 \}$, 
    \quad 
    $\mathcal N_3 = \textrm{span} \{ (0,z_2), (0,z_1), B_0 \}$;  \\ 
    $\mathcal N_p = \mathcal N_3 \ \textrm{for any} \ p \geqslant 3$. 
   \end{tabular} 
  \end{center} 
  The first three follow from $Qz=C$, $Qz_1=-z$, $Qz_2=-z_1$, and 
  $C \cdot z = 0$, $C \cdot z_1 = 1$, $C \cdot z_2 = 0$.  The fourth is 
  essentially due to the fact that, while certainly there is a $z_3$ such 
  that $Qz_3=-z_2$, it is unable to satisfy $C \cdot z_3 = 0$ because of (*) 
  above.  The union of all the $\mathcal N_i$ is the generalized null space 
  $\mathcal N$ of $A$.  It is invariant under $A$. 

  Normalize $(0,z_1)$, $(0,z_2)$ to yield two real $\langle \, , \, \rangle$ 
  orthonormal spacelike vectors $X_1$, $X_2$.  A routine calculation 
  produces the unit timelike real vector
  \begin{displaymath} 
   X_0 := \frac{|z_2|}{|z \cdot z_2|} (1,z) + X_2 \, ,  
  \end{displaymath} 
  which is $\langle \, , \, \rangle$ orthogonal to $X_1$, $X_2$.  Also,  
  with $a_1 := |z_1| / |z_2|$, we have $A X_0 = ia_1 X_1$, 
  $A X_1 = ia_1 (X_0-X_2)$, and $A X_2 = ia_1 X_1$.  Let $B_1$ be the 
  real $\langle \, , \, \rangle$ orthonormal basis 
  $\{ X_0, X_1, X_2, (0,e_j), j = 2h+1, \dots, n \}$ for the generalized 
  null space $\mathcal N$.  With respect to $B_1$, the matrix of 
  $A_{|\mathcal N}$ has the form $ia_1 T \oplus 0_{n-2h}$, where 
  \begin{displaymath} 
   T = 
   \left(
    \begin{array}{ccc}
     0  &  1  &  0  \\
     1  &  0  &  1  \\
     0  & -1  &  0  \\
    \end{array}
   \right) \, . 
  \end{displaymath} 
  Correspondingly, the matrix of $\Omega_{|\mathcal N}$  
  is $a_1 T \oplus 0_{n-2h}$, with $a_1 > 0$.  

  Since $X_0$ is timelike, the scalar product becomes positive definite on 
  $\mathcal N^\perp$, which is invariant under the self-adjoint $A$.  
  By \S 10.2, there is a real orthonormal basis $B_2$ for $\mathcal N^\perp$ 
  which puts $A_{|\mathcal N^\perp}$, and hence $\Omega_{|\mathcal N^\perp}$, 
  into normal form.  Incidentally, this normal form must look like 
  $a_2 J \oplus \cdots \oplus a_{m'} J$, 
  where $a_2 \geqslant \cdots \geqslant a_{m'} > 0$, because the kernel of 
  $\Omega$ has already been accounted for in $\mathcal N$. 

  Since $X_0$ is future-pointing, the real $\langle \, , \, \rangle$ 
  orthonormal basis $\mathcal B := B_1 \cup B_2$ gives an element 
  $g := {\mathcal B}^{-1} \in O_+(1,n)$.  The normal form $\tilde \Omega 
  := g \Omega g^{-1}$ is as follows:  
    \begin{center} 
   \begin{tabular}{lll} 
    when $n$ is even, 
    & 
    $\tilde \Omega = a_1 T \oplus a_2 J \oplus \cdots \oplus a_m J$, 
    & 
    with $m = n / 2$; \\ 
    when $n$ is odd, 
    &
    $\tilde \Omega = a_1 T \oplus a_2 J \oplus \cdots \oplus a_m J \oplus 0$, 
    &
    with $m = (n-1)/2$.  
   \end{tabular}  
  \end{center} 
  Here, $a_1 > 0$ and $a_2 \geqslant \cdots \geqslant a_m \geqslant 0$. 

%------------------------------------------------------------------------------

%------------------------------------------------------------------------------

\end{document}